\documentclass[11pt]{amsart}
\usepackage{lmodern}
\usepackage{amsmath,amssymb,amsthm,mathtools}
\usepackage{microtype}
\usepackage{enumitem}
\usepackage{aliascnt}
\usepackage{hyperref}
\usepackage[nameinlink,capitalise,noabbrev]{cleveref}
\usepackage[margin=1.12in]{geometry}
\usepackage{tikz-cd}
\hypersetup{colorlinks=true,linkcolor=blue,citecolor=blue,urlcolor=blue}

\newcommand{\Z}{\mathbb Z}

\newcommand{\C}{\mathbb C}
\newcommand{\R}{\mathbb R}

\newcommand{\CP}{\mathbb {CP}}
\newcommand{\eps}{\varepsilon}
\newcommand{\id}{\operatorname{id}}

\newcommand{\im}{\operatorname{im}}

\newcommand{\PL}{\mathrm{PL}}

\newcommand{\BO}{\mathrm{BO}}
\newcommand{\BPL}{\mathrm{BPL}}

\newcommand{\MO}{\mathrm{MO}}
\newcommand{\MPL}{\mathrm{MPL}}

\newcommand{\Diff}{\mathrm{Diff}}

\newcommand{\Ric}{\mathrm{Ric}}
\newcommand{\Vol}{\mathrm{Vol}}

\theoremstyle{plain}
\newtheorem{theorem}{Theorem}[section]
\newaliascnt{proposition}{theorem}
\newtheorem{proposition}[proposition]{Proposition}
\aliascntresetthe{proposition}
\newaliascnt{lemma}{theorem}
\newtheorem{lemma}[lemma]{Lemma}
\aliascntresetthe{lemma}
\newaliascnt{corollary}{theorem}
\newtheorem{corollary}[corollary]{Corollary}
\aliascntresetthe{corollary}
\theoremstyle{definition}
\newaliascnt{definition}{theorem}

\aliascntresetthe{definition}
\theoremstyle{remark}
\newaliascnt{remark}{theorem}
\newtheorem{remark}[remark]{Remark}
\aliascntresetthe{remark}

\crefname{theorem}{Theorem}{Theorems}
\crefname{proposition}{Proposition}{Propositions}
\crefname{lemma}{Lemma}{Lemmas}
\crefname{corollary}{Corollary}{Corollaries}
\crefname{definition}{Definition}{Definitions}
\crefname{remark}{Remark}{Remarks}

\title[Linear circle quotients of products of odd spheres]
{Linear circle quotients of products of odd-dimensional spheres: classification and geometric applications}

\author[Wen Shen]{Wen Shen}
\address{Department of Mathematics, Wenzhou University,
Wenzhou, China}
\email{shenwen121212@163.com}

\author[Yang Su]{Yang Su}
\address{School of Mathematics and Systems Science, Chinese Academy of Sciences, Beijing, China; \newline
               School of Mathematical Sciences, University of Chinese Academy of Sciences, Beijing, China}
 \email{suyang@math.ac.cn}

\subjclass[2020]{Primary 57R15, 57R19, 57R55, 57R65; Secondary 55N22, 55R25, 55Q45, 53C25, 53D10, 53E20}
\keywords{circle action, orbit manifold, modified surgery, normal invariant, Pontryagin class, real K-theory, sphere bundle, homogeneous Einstein metric, contact structure}

\begin{document}

\begin{abstract}
In this paper, we study the orbit manifolds of free diagonal linear circle actions on \(S^{2n+1}\times S^{2n+1}\). For every \(n\geq 1\), we establish a  PL homeomorphism classification of the orbit manifolds in terms of the weights of the actions. We also obtain diffeomorphism classification and homotopy classification of such orbit manifolds in dimension $9$ and $13$.  The diffeomorphism classification results lead to two geometric applications. For every \(m\geq 1\), there exist a non-spin \(9\)-manifold \(M_m^9\) and a spin \(13\)-manifold \(M_m^{13}\) in our family such that the moduli space of unit-volume homogeneous Einstein metrics on each manifold has at least \(2^{m-1}\) connected components. Moreover, each of these manifolds admits at least \(2^{m-1}\) regular Sasakian structures whose underlying contact structures are pairwise non-contactomorphic.
 \end{abstract}

\maketitle

\section{Introduction}

Circle quotients and biquotients have played
an important role in the construction and classification of
manifolds with positive sectional curvature.  Classical examples
include the Aloff--Wallach spaces and Eschenburg spaces, which
are parametrized by integral weights
\cite{AloffWallach1975,Eschenburg1982}.  The classification of
Aloff--Wallach spaces by Kreck and Stolz demonstrates how arithmetic information
in the defining parameters can distinguish manifolds with
different smooth structures
\cite{KreckStolz1991,KreckStolz1998}.  In dimension thirteen,
the Bazaikin spaces provide another important family of
curvature-related biquotients, whose topological invariants
can also be expressed explicitly in terms of their defining
integer parameters \cite{FloritZiller2009}.

  Another particularly relevant example is the construction
of homogeneous Einstein metrics on principal torus bundles by
Wang and Ziller \cite{WangZiller1990}.  In particular, their
construction yields homogeneous Einstein metrics on certain
principal circle bundles over products of compact homogeneous
K\"ahler--Einstein manifolds.  The case of circle bundles over
$\CP^n\times\CP^n$ is especially relevant to the present work:
their total spaces can be realized as quotients of
$S^{2n+1}\times S^{2n+1}$ by diagonal linear $S^1$-actions
with repeated integral weights.  Thus, the integer parameters
defining the circle action also determine a family of manifolds
admitting homogeneous Einstein metrics.  This provides a direct
connection between the arithmetic of circle actions and the
geometry of their orbit manifolds.

In this paper, we investigate the orbit manifolds of general free linear circle actions on a product of spheres $S^{2n+1} \times S^{2n+1}$.  To be precise, let
\[
 a=(a_1,\ldots,a_{n+1}),\qquad
 b=(b_1,\ldots,b_{n+1})\in\mathbb Z^{n+1},
\]
and consider the diagonal linear action of \(S^1\) on \(S^{2n+1}\times S^{2n+1}\),  
\begin{equation}\label{eq:action}
z\cdot((u_1,\ldots,u_{n+1}),(v_1,\ldots,v_{n+1}))
=((z^{a_1}u_1,\ldots,z^{a_{n+1}}u_{n+1}),
(z^{b_1}v_1,\ldots,z^{b_{n+1}}v_{n+1})).
\end{equation}
 The action is free if and only if (see \cref{prop:free})
\[
 \gcd(|a_i|,|b_j|)=1\qquad\text{for all }i,j.
\]
Whenever the action is free,
we denote its orbit manifold by
\[
M_{a,b}
=
\bigl(S^{2n+1}\times S^{2n+1}\bigr)/S^1.
\]

\begin{remark}
It is enlightening to contrast this situation with actions on a single odd-dimensional
sphere.  For a diagonal linear action of $S^1$ on $S^{2n+1}$,
\[
 z\cdot(w_1,\ldots,w_{n+1})
 =(z^{c_1}w_1,\ldots,z^{c_{n+1}}w_{n+1}),
\]
freeness forces \(c_i= \pm 1\) for every \(i\).  After
complex conjugating the coordinates with negative weight, such an action is
equivalent to the standard Hopf action and its orbit space is \(\CP^n\).  If
finite stabilizers are allowed, general weights lead instead to weighted
projective spaces; see \cite{Kawasaki1973}.  On a product of two spheres the
freeness condition is much less rigid. Thus products of spheres provide the first natural setting in which free
linear circle actions retain genuinely nontrivial arithmetic parameters.
\end{remark}

\begin{remark} Ozeki and Uchida in 
\cite{OzekiUchida1972} investigated principal circle
actions on products of spheres and computed the integral
cohomology rings and Pontryagin classes of the corresponding
orbit manifolds.  While these invariants provide information
about the orbit manifolds, our main objective is to understand
their classification up to homotopy equivalence, PL homeomorphism or diffeomorphism.
\end{remark}

\begin{remark}
The action in (\ref{eq:action}) is isometric for the product round metric, so the
O'Neill formula gives every \(M_{a,b}\) a metric of nonnegative sectional
curvature \cite{ONeill1966}.  The classification therefore takes place
inside a geometrically distinguished class of manifolds.
\end{remark}

Now we introduce the numerical topological invariants appearing in the classification theorems, which are determined by the weights $(a,b)$. For details see  \cref{sec:orbit}. First noticte that all the quotients have the same integral cohomology ring,
\[
H^*(M_{a,b};\mathbb Z)
\cong
H^*(\CP^n\times S^{2n+1};\mathbb Z).
\]
Thus ordinary cohomology does not distinguish the members of the family.
The weights nevertheless determine the stable tangent bundle and hence the
characteristic and normal data of the quotient.  To be more precise, for
\(1\le k\le\lfloor n/2\rfloor\), put
\begin{equation}\label{eq:Pi}
 \Pi_k(a,b)
 =
 e_k(a_1^2,\ldots,a_{n+1}^2,b_1^2,\ldots,b_{n+1}^2),
\end{equation}
where \(e_k\) denotes the \(k\)th elementary symmetric polynomial.  Then
\[
 p_k(M_{a,b})=\Pi_k(a,b)x^{2k}.
\]
where $x$ is a generator of $H^2(M_{a,b};\mathbb Z)$.
When \(n=4r+1\), there is an additional invariant of the stable tangent bundle. Namely, for \(m\ge0\), set
\begin{equation}\label{eq:cjm}
 c_j(m)=
 \begin{cases}
 \displaystyle \frac{m}{j}\binom{m+j-1}{2j-1},&m\ge j,\\[1.2ex]
 0,&m<j,
 \end{cases}
\end{equation}
and define
\begin{equation}\label{eq:Delta}
 \Delta_r(a,b)=
 \sum_{i=1}^{n+1}c_{2r+1}(|a_i|)
 +\sum_{j=1}^{n+1}c_{2r+1}(|b_j|)
 \pmod2.
\end{equation}

\begin{remark}
Here we briefly explain the meaning of the invariant $\Delta_r(a,b)$. For details see section \ref{subsec:projective-stage}. The classifying map
$
f:M_{a,b}\longrightarrow\CP^\infty
$
of the principal circle bundle
\(S^{2n+1}\times S^{2n+1}\to M_{a,b}\) is
\((2n+1)\)-connected.  Consequently there is a canonical projective stage
\begin{equation}\label{eq:stage}
i:\CP^n\longrightarrow M_{a,b},
\end{equation}
unique up to the sign of the degree-two generator. When $n=4r+1$, the Pontryagin classes do not contain all of the
stable real bundle information on the projective stage. The invariant $\Delta_{a,b}$ is exactly the top torsion coefficient in Fujii's normal
form of \(KO(\CP^n)\).  It is the only stable real bundle information in this family that is
invisible to the Pontryagin classes.
\end{remark}
 
Our main topological classification is the following.

\begin{theorem}\label{thm:main}
Let \(n\ge1\), and let \(M_{a,b}\) and \(M_{a',b'}\) be orbit manifolds of
free diagonal linear \(S^1\)-actions on
\(S^{2n+1}\times S^{2n+1}\).  Then $M_{a,b}$ is \(\PL\) homeomorphic to $M_{a',b'}$
if and only if
\[
 \Pi_k(a,b)=\Pi_k(a',b')
 \qquad
 (1\le k\le\lfloor n/2\rfloor),
\]
and, when \(n=4r+1\), also
\[
 \Delta_r(a,b)=\Delta_r(a',b').
\]
\end{theorem}


\begin{remark}
When \(n=1\), the list of Pontryagin invariants is empty and
\[
 \Delta_0(a,b)=a_1^2+a_2^2+b_1^2+b_2^2\pmod2.
\]
The five-dimensional case therefore recovers the spin/non-spin distinction:
\(M_{a,b}\) is diffeomorphic to \(S^2\times S^3\) when \(\Delta_0=0\), and
to the nontrivial oriented \(S^3\)-bundle
\(S^2\widetilde\times S^3\) when \(\Delta_0=1\).  

When \(n=3\), Lagrange's four-square theorem implies that  
every \(13\)-dimensional quotient is PL homeomorphic to a linear quotient
\(M_{\mathbf1,q}\), where \(\mathbf1=(1,1,1,1)\), $q=(q_1,q_2,q_3,q_4)$.  
\end{remark}

\medskip

A consequence of the PL classification theorem is a sphere bundle realization which gives an alternative description of the quotient
manifolds and provides a passage from the PL homeomorphism classification to homotopy classification.

\begin{theorem}\label{thm:sphere-intro}
Let \(n\ge2\).  Every orbit manifold \(M_{a,b}\) is $\PL$ homeomorphic to the
unit sphere bundle of an oriented real vector bundle of rank \(2n+2\) over
\(\CP^n\).
\end{theorem}

Before continuing to present the smooth classification results, we pause to ponder the framework demonstrated by  Theorem \ref{thm:main} and \ref{thm:sphere-intro}. The projective stage (\ref{eq:stage}) separates the
two main levels of the classification.  Let \(L\to\CP^\infty\) be the
universal complex line bundle, let \(j_n:\CP^n\to\CP^\infty\) be the
canonical inclusion, and put
\[
 \gamma_{a,b}
 =
 \left[
 j_n^*\left(
 \bigoplus_{i=1}^{n+1}L^{a_i}
 \oplus
 \bigoplus_{j=1}^{n+1}L^{b_j}
 \right)_{\!\mathbb R}
 \right]
 -[\varepsilon^{4n+4}]
 \in\widetilde{KO}(\CP^n).
\]
The stable tangent bundle formula of $M_{a,b}$ (cf. Proposition \ref{prop:tangent}) identifies \(\gamma_{a,b}\) with the
reduced stable tangent class on the projective stage.  Conceptually, the PL
classification Theorem \ref{thm:main} may be summarized by
\begin{equation}\label{eq:intro-KO-PL}
 M_{a,b}\cong_{\PL}M_{a',b'}
 \quad\Longleftrightarrow\quad
 \gamma_{a,b}=\gamma_{a',b'}
 \quad\text{in }\widetilde{KO}(\CP^n).
\end{equation}
On the other hand, by the sphere bundle realization Theorem \ref{thm:sphere-intro}, passage from PL classification to homotopy classification is governed by the stable \(J\)-homomorphism, which measures the stable fiber homotopy class of this sphere bundle, leading to the criterion
\begin{equation}\label{eq:intro-J-homotopy}
 M_{a,b}\simeq M_{a',b'}
 \quad\Longleftrightarrow\quad
 J(\gamma_{a,b})=J(\gamma_{a',b'})
 \quad\text{in }\widetilde{J}(\CP^n)
\end{equation}
where $\widetilde{J}(\CP^{n})$ denotes the classical reduced $J$-group \cite{Adams1963} over $\CP^{n}$, the $J$-homomorphism assigns to a
stable real vector bundle its associated stable
spherical fibration.
In principle the
criterion can be converted into finitely many congruences in the weights by
combining the Chebyshev expansion \eqref{eq:cjm} with the relations in
\(\widetilde{J}(\CP^{n})\).  In dimensions \(9\) and \(13\), 
\cite[Lemma 1.4]{Adams1966} reduces the criterion to the single congruence modulo
\(24\) appearing in \cref{thm:9d-intro,thm:13d-intro}.

The basic tool of the PL classification is Kreck's modified surgery theory, which upgrades equality of the stable real bundle data to PL
homeomorphism without first solving the homotopy classification problem --- in our situation the more complicated
\(\widetilde{J}(\CP^{n})\)-problem. The PL classification is governed by a single normal bordism framework.  For
\(n\ge2\), the weight data determine the stable normal bundle on the
projective stage and hence a common normal \(2n\)-type.  The equivariant
coordinate projections of the linear action force the canonical normal
smoothing into a low Atiyah--Hirzebruch filtration, and this filtration
disappears after passage to PL bordism.  By modified surgery theory one then obtains a PL $h$-cobordism between the orbit manifolds.  

\medskip 

The smooth classification of the orbit manifolds is subtler. For small values $n=2$ and $3$, we obtain an almost complete classification. For \(n=2,3\), write
\[
\lambda(a,b)=\sum_i(a_i^2+b_i^2)=\Pi_1(a,b).
\]

\begin{theorem}\label{thm:9d-intro}
Let \(M_{a,b}\) and \(M_{a',b'}\) be orbit manifolds of free diagonal
linear \(S^1\)-actions on \(S^5\times S^5\). Then
\begin{enumerate}[label=\textup{(\roman*)}]
\item
$ M_{a,b}\simeq M_{a',b'}
 \quad\Longleftrightarrow\quad
 \lambda(a,b)\equiv\lambda(a',b')\pmod{24}.$
\item If both quotients are non-spin, equivalently if both
\(\lambda\)-invariants are odd, then
\[
 M_{a,b}\cong_{\Diff}M_{a',b'}
 \quad\Longleftrightarrow\quad
 \lambda(a,b)=\lambda(a',b').
\]
\end{enumerate}
\end{theorem}

\begin{theorem}\label{thm:13d-intro}
Let \(M_{a,b}\) and \(M_{a',b'}\) be orbit manifolds of free diagonal
linear \(S^1\)-actions on \(S^7\times S^7\). Then
\begin{enumerate}[label=\textup{(\roman*)}]
\item
$ M_{a,b}\simeq M_{a',b'}
 \quad\Longleftrightarrow\quad
 \lambda(a,b)\equiv\lambda(a',b')\pmod{24}.$
\item If both quotients are non-spin, equivalently if both
\(\lambda\)-invariants are odd, or if both $\lambda$-invariants are divisible by
\(4\), then
\[
 M_{a,b}\cong_{\Diff}M_{a',b'}
 \quad\Longleftrightarrow\quad
 \lambda(a,b)=\lambda(a',b').
\]
\end{enumerate}
\end{theorem}

In
dimension \(9\), the proof is based on an analysis of the smooth surgery exact sequence of the orbit manifold. The non-spin condition forces a stable \(\eta\)-component
in the top-cell attaching map in its CW-structure, which kills the smooth normal invariant. The
remaining \(L_{10}(e)\)-ambiguity is removed by a calculation with the surgery obstruction formula of Rourke-Sullivan.  In dimension \(13\), the non-spin
case is proved by modified surgery theory, with a detailed computation of the corresponding string bordism group; whereas
the \(\lambda\equiv0\pmod4\) case is proved by the smooth surgery exact sequence, based on the PL classification. The
remaining spin congruence class \(\lambda\equiv2\pmod4\) is not pursued here.

The smooth classification also produces geometric multiplicity on fixed smooth
manifolds.  Consider the repeated-weight subfamily
\[
M_{k,l}^{(n)}
=
M_{(k,\ldots,k),(l,\ldots,l)}.
\]
For coprime nonzero \(k,l\), this is a principal circle bundle over
\(\CP^n\times\CP^n\).  Wang and Ziller constructed homogeneous Einstein
metrics on these bundles \cite{WangZiller1990}.  When \(k,l>0\), the Euler
class is represented by an integral K\"ahler form, so the Boothby--Wang
construction also gives a regular Sasakian structure
\cite{BoothbyWang1958,BoyerGalicki2008}.  Our smooth classification allows
different circle-bundle presentations to be transported to the same
smooth manifold.

Let \(\mathcal E_h(M)\) denote the moduli space of unit-volume homogeneous
Einstein metrics on \(M\), with no transitive group
fixed in advance.

\begin{theorem}[Geometric multiplicity]\label{thm:Einstein-intro}
For every \(m\ge1\), there exist a non-spin \(9\)-dimensional orbit manifold
\(M_m^9\) and a spin \(13\)-dimensional orbit manifold \(M_m^{13}\) such
that
\[
 \#\pi_0\mathcal E_h(M_m^9)\ge2^{m-1},
 \qquad
 \#\pi_0\mathcal E_h(M_m^{13})\ge2^{m-1}.
\]
Moreover, each of \(M_m^9\) and \(M_m^{13}\) carries at least
\(2^{m-1}\) regular Sasakian structures whose underlying cooriented contact
structures are pairwise non-contactomorphic.
\end{theorem}

The source of the geometric multiplicity is an elementary number-theoretic multiplicity.  If
\[
 N_m=p_1\cdots p_m,
\]
with \(p_i\) distinct primes and $p_i\equiv 1 \pmod4$, then \(N_m\) has exactly \(2^{m-1}\) primitive
 representations 
\[
N_m=k^2+l^2 \ \ \ \  \textrm{with \(0<k<l\).}
\]  In dimensions
\(9\) and \(13\), the smooth classification identifies all corresponding
quotients with one manifold.  The normalized Wang--Ziller Einstein constants
separate the resulting Einstein metrics, while the first Chern classes of
the Boothby--Wang contact distributions separate the regular Sasakian
structures.  

\medskip

The paper is organized as follows.  In \cref{sec:orbit} we determine the
freeness condition, compute the cohomology and stable tangent bundle, and
identify the additional \(KO\)-torsion invariant in
\cref{subsec:projective-stage}.  In \cref{sec:normal-type} we construct a
common normal \(2n\)-type for quotients with the same stable bundle data.
\Cref{sec:smooth-filtration} studies the corresponding smooth normal bordism
classes, with the equivariant prime-avoidance factorization isolated in
\cref{sec:prime}.  We pass to the PL normal type in \cref{sec:PL} and complete
the PL classification in \cref{sec:surgery} using modified surgery theory.
The low-dimensional smooth refinements occupy the next two sections.
\Cref{sec:13smooth} treats dimension \(13\), including the non-spin case and
the spin class \(\lambda\equiv0\pmod4\); \cref{sec:9smooth} treats the
non-spin \(9\)-dimensional case.  \Cref{sec:sphere-realization} proves the
sphere bundle realization and the stable \(J\)-theoretic homotopy criterion,
with explicit consequences in dimensions \(9\) and \(13\).  Finally,
\cref{sec:Einstein} develops the applications to homogeneous Einstein and Sasakian structures.
.

\section{Orbit manifolds and characteristic classes}\label{sec:orbit}

Let $P=S^{2n+1}\times S^{2n+1}$ and suppose that \eqref{eq:action} is free.
Write $p:P\to M_{a,b}$ for the quotient map.

\begin{proposition}\label{prop:free}
The action \eqref{eq:action} is free if and only if
\[
 \gcd(|a_i|,|b_j|)=1
 \qquad\text{for every }1\le i,j\le n+1,
\]
where $\gcd(0,m)=|m|$.
\end{proposition}

\begin{proof}
At a point supported on the $i$th coordinate of the first sphere and the
$j$th coordinate of the second, the stabilizer is
\[
 \{z\in S^1:z^{a_i}=z^{b_j}=1\}.
\]
This group is trivial precisely when $\gcd(|a_i|,|b_j|)=1$.
If $a_i=b_j=0$, it is the entire circle; we use $\gcd(0,0)=0$.
At an arbitrary point, choose one nonzero coordinate from each sphere.
Its stabilizer is contained in the corresponding group above.  These
observations prove both implications.
\end{proof}

Put
\begin{equation}\label{eq:AB}
 A(a)=\prod_{i=1}^{n+1}a_i,
 \qquad
 B(b)=\prod_{j=1}^{n+1}b_j.
\end{equation}
By \cref{prop:free},
\begin{equation}\label{eq:gcdAB}
 \gcd(|A(a)|,|B(b)|)=1.
\end{equation}
In particular, at least one of $A(a),B(b)$ is nonzero.

Let $x\in H^2(\CP^\infty;\Z)$ be the canonical generator.  The Borel
construction of the free action is homotopy equivalent to the quotient:
\[
 ES^1\times_{S^1}P\simeq M_{a,b}.
\]
The associated fibration
\begin{equation}\label{eq:borel}
 P\longrightarrow M_{a,b}\xrightarrow{f}\CP^\infty
\end{equation}
has Serre spectral sequence
\[
 E_2^{p,q}=H^p(\CP^\infty;H^q(P;\Z)).
\]
Let $u,v\in H^{2n+1}(P;\Z)$ be the fundamental classes of the two sphere
factors.

\begin{lemma}\label{lem:transgression}
In the spectral sequence of \eqref{eq:borel},
\[
 d_{2n+2}(u)=A(a)x^{n+1},
 \qquad
 d_{2n+2}(v)=B(b)x^{n+1}.
\]
\end{lemma}

\begin{proof}
The projection onto the first factor induces a map from \eqref{eq:borel} to
the sphere bundle over $\CP^\infty$ associated with the complex representation
\[
 \C_{a_1}\oplus\cdots\oplus\C_{a_{n+1}}.
\]
The transgression of the sphere fundamental class is the Euler class, which is
the top Chern class
$A(a)x^{n+1}$.  The second formula follows in the same way.
\end{proof}

\begin{corollary}\label{cor:cohomology}
The manifold $M_{a,b}$ is simply connected and
\[
 H^*(M_{a,b};\Z)\cong
 \Z[x,y]/(x^{n+1},y^2),
 \qquad |x|=2,\quad |y|=2n+1.
\]
Moreover, the map $f$ in \eqref{eq:borel} is $(2n+1)$-connected.
\end{corollary}

\begin{proof}
By \eqref{eq:gcdAB}, the homomorphism
\[
 \Z\{u,v\}\longrightarrow\Z\{x^{n+1}\},
 \qquad (r,s)\longmapsto (rA(a)+sB(b))x^{n+1},
\]
is surjective.  Its kernel has rank one.  The multiplicative spectral sequence
therefore gives the stated ring, with an odd generator represented by a
primitive generator of that kernel.  The homotopy exact sequence of the principal circle bundle
\(P\to M_{a,b}\) shows that \(M_{a,b}\) is simply connected.
The homotopy fiber of $f$ is $P$, which is $2n$-connected.
The homotopy exact sequence therefore shows directly that $f$ is
$(2n+1)$-connected.
\end{proof}

We next record an equivariant tangent-bundle identity.  Let $L\to\CP^\infty$
be the universal complex line bundle, with the associated-bundle convention
$P\times_{S^1}\C_m\cong f^*L^m$, and set
\begin{equation}\label{eq:Gamma}
 \Gamma_{a,b}=\bigoplus_{i=1}^{n+1}L^{a_i}
 \oplus\bigoplus_{j=1}^{n+1}L^{b_j}.
\end{equation}
Let $\Gamma_{a,b,\R}$ denote its underlying real bundle.

\begin{proposition}\label{prop:tangent}
There is a stable bundle isomorphism over $M_{a,b}$
\begin{equation}\label{eq:tangent}
 TM_{a,b}\oplus\eps^3\cong f^*\Gamma_{a,b,\R}.
\end{equation}
Consequently,
\begin{equation}\label{eq:pontryagin-total}
 p(TM_{a,b})=
 \prod_{i=1}^{n+1}(1+a_i^2x^2)
 \prod_{j=1}^{n+1}(1+b_j^2x^2),
\end{equation}
where $x=f^*(x)\in H^2(M_{a,b};\Z)$.
\end{proposition}

\begin{proof}
Write
$$
P=S(V_a)\times S(V_b),
\qquad
V_a=\bigoplus_{i=1}^{n+1}\C_{a_i},
\qquad
V_b=\bigoplus_{j=1}^{n+1}\C_{b_j},
$$
where \(\C_m\) denotes the one-dimensional complex representation on
which \(z\in S^1\) acts by multiplication by \(z^m\).

For a unitary \(S^1\)-representation \(V\), the unit sphere
\(S(V)\subset V_{\R}\) has the standard radial splitting
$$
TS(V)\oplus\eps^1
\longrightarrow
S(V)\times V_{\R},
\qquad
(u;\xi,t)\longmapsto (u,\xi+tu).
$$
This is an \(S^1\)-equivariant bundle isomorphism: indeed, for
\(z\in S^1\),
$$
z(\xi+tu)=z\xi+t(zu),
$$
and the unitary action sends \(T_uS(V)\) isomorphically onto
\(T_{zu}S(V)\).  Applying this to \(V_a\) and \(V_b\) and taking the
product gives an \(S^1\)-equivariant isomorphism
\begin{equation}\label{eq:equivariant-radial}
TP\oplus\eps^2
\cong
P\times(V_a\oplus V_b)_{\R}.
\end{equation}

Next consider the principal circle bundle
$$
S^1\longrightarrow P\xrightarrow{p}M_{a,b}.
$$
The fundamental vector field
\[
 Z(q)=\left.\frac{d}{dt}\right|_{t=0}e^{it}\cdot q \quad (q\in P, \quad t\in \mathbb{R})
\]
is nowhere zero because the action is free.  Since $S^1$ is abelian,
$d(A_z)_qZ(q)=Z(A_zq)$, where $A_z$ is the action of $z$ on $P$.
Thus $\ker(dp)$ is equivariantly trivial.  The product round metric is invariant, and its
orthogonal horizontal distribution identifies equivariantly with $p^*TM$.
Consequently,
\begin{equation}\label{eq:vertical-splitting}
TP\cong p^*TM_{a,b}\oplus\eps^1
\end{equation}
is an $S^1$-equivariant isomorphism.
Combining \eqref{eq:equivariant-radial} and
\eqref{eq:vertical-splitting}, we obtain
$$
p^*TM_{a,b}\oplus\eps^3
\cong
P\times(V_a\oplus V_b)_{\R}
$$
as \(S^1\)-equivariant real vector bundles over \(P\).

Passing to the quotient by the free circle action gives
$$
TM_{a,b}\oplus\eps^3
\cong
P\times_{S^1}(V_a\oplus V_b)_{\R}.
$$

Since \(f:M_{a,b}\to BS^1=\CP^\infty\) classifies the principal bundle
\(P\to M_{a,b}\), we have
$$
P\times_{S^1}(V_a\oplus V_b)_{\R}
\cong
f^*
\left(
\bigoplus_{i=1}^{n+1}L^{a_i}
\oplus
\bigoplus_{j=1}^{n+1}L^{b_j}
\right)_{\R}
=
f^*\Gamma_{a,b,\R},
$$
which proves \eqref{eq:tangent}.
Finally, the Whitney product formula gives \eqref{eq:pontryagin-total}.
\end{proof}

In particular,
\begin{equation}\label{eq:pk}
 p_k(M_{a,b})=\Pi_k(a,b)x^{2k},
 \qquad 1\le k\le\lfloor n/2\rfloor.
\end{equation}

\begin{remark}\label{rem:n1-w2}
The stable tangent formula above is also valid for $n=1$.  If $\bar x$ is
the mod-$2$ reduction of the generator of $H^2(M_{a,b};\mathbb Z)$, then
\[
 w_2(M_{a,b})=(a_1+a_2+b_1+b_2)\bar x
             =\Delta_0(a,b)\bar x.
\]
Thus the $n=1$ invariant in \cref{thm:main} is exactly the spin/non-spin
invariant in the Smale--Barden classification.
\end{remark}

\subsection{\texorpdfstring{The projective stage and the $KO$ torsion coefficient}{The projective stage and the KO torsion coefficient}}\label{subsec:projective-stage}

By \cref{cor:cohomology}, after replacing $f$ by a fibration, obstruction
theory gives a map
\begin{equation}\label{eq:projective-stage}
 i:\CP^n\longrightarrow M_{a,b}
 \qquad\text{such that}\qquad
 f\circ i\simeq j_n:\CP^n\longrightarrow\CP^\infty,
\end{equation}
where $j_n$ is the canonical inclusion.  Since the homotopy fiber of $f$ is
$2n$-connected, the lift $i$ is unique up to homotopy after the sign of $x$ is
fixed.

We use the following consequence of the calculations of Fujii and Sanderson.

\begin{lemma}\label{lem:KOCP}
Let $n\ge2$.
\begin{enumerate}[label=(\roman*)]
\item The restriction $KO(\CP^\infty)\to KO(\CP^n)$ is surjective.
\item If $n\not\equiv1\pmod4$, a stable real bundle of fixed rank over $\CP^n$ is determined
by its Pontryagin classes.
\item If \(n=4r+1\), let
\[
r_{\mathbb R}:K(\CP^n)\longrightarrow KO(\CP^n)
\]
denote realification and put
\(\mu=r_{\mathbb R}(L-1)\in\widetilde{KO}(\CP^n)\).  Then
\[
 \widetilde{KO}(\CP^{4r+1})
 \cong \Z\{\mu,\ldots,\mu^{2r}\}
 \oplus (\Z/2)\{\mu^{2r+1}\}.
\]
Thus, in this congruence class a stable real bundle is determined by its
Pontryagin classes and the coefficient of $\mu^{2r+1}$.
\end{enumerate}
\end{lemma}

\begin{proof}
This is the standard form of Fujii's computation together with Sanderson's
stable bundle classification; see \cite{Fujii1967,Sanderson1964}.  The reduced group is generated additively by powers of
$\mu=[L_{\R}]-2$ (with the stated relation in the torsion case).
Each such polynomial extends over $\CP^\infty$, which proves the
surjectivity assertion.  All stable bundle comparisons here have fixed rank,
or equivalently are comparisons of reduced classes.
\end{proof}

When $n=4r+1$, let $\delta(\zeta)\in\Z/2$ denote the coefficient of
$\mu^{2r+1}$ in Fujii's normal form for a stable real bundle $\zeta$ over
$\CP^n$, and define
\begin{equation}\label{eq:deltaM}
 \delta(M_{a,b})=\delta(i^*TM_{a,b}).
\end{equation}

\begin{proposition}\label{prop:delta-intrinsic}
Let $n=4r+1$.  The coefficient
\[
 \delta(M)=\delta(i_M^*TM)\in\Z/2
\]
is independent of the choices in the projective-stage construction and is
preserved by homeomorphisms.  More generally, it is preserved by a homotopy
equivalence that identifies the rational Pontryagin classes.  
\end{proposition}

\begin{proof}
Put $K=\widetilde{KO}(\CP^{4r+1})$ and
$T=\operatorname{Tor}K=(\Z/2)\{\mu^{2r+1}\}$.
We first verify the precise $J$-theoretic fact that is needed:
\begin{equation}\label{eq:J-injective-torsion}
 \ker(J:K\longrightarrow\widetilde{J}(\CP^{4r+1}))\cap T=0.
\end{equation}

Indeed, suppose that the nonzero element
\(\tau=\mu^{2r+1}\in T\) lies in \(\ker J\).  Since \(\tau\) is
\(2\)-primary, Obiedat's localization formula
\cite[Corollary~2.4]{Obiedat1999Localization} gives
$$
 \tau=(\psi^3-1)v
 \qquad\text{for some }v\in K_{(2)}.
$$
Passing to the torsion-free quotient \(K/T\), the operator
\(\psi^3-1\) is injective: rationally, the Pontryagin character
identifies its degree-\(4j\) eigenspace with eigenvalue \(3^{2j}-1\),
which is nonzero for every \(j\ge1\).  Hence \(v\in T\).

On the other hand, since
$\psi^3(\mu)
=9\mu +6\mu^2+\mu^3$ and $\mu^{2r+2}=0$, we have
$$
 \psi^3(\mu^{2r+1})
 =
 3^{4r+2}\mu^{2r+1}
 =
 \mu^{2r+1}.
$$
Here note that \(T\) has order two.  Thus
$$
 (\psi^3-1)v=0,
$$
contradicting \((\psi^3-1)v=\tau\ne0\).  This proves
\eqref{eq:J-injective-torsion}.

For a fixed degree-two generator, the projective-stage lift is unique
up to homotopy.  Changing its sign precomposes the lift with complex
conjugation, which fixes $\mu$ and all of its powers.  Thus $\delta(M)$
is well defined.

Now let $h:M\to M'$ be a homotopy equivalence identifying the rational
Pontryagin classes.  After adjusting the sign of a generator, the projective
stages satisfy $h\circ i_M\simeq i_{M'}$.  Invariance of the Spivak normal
spherical fibration gives
\[
 J(i_M^*TM)=J(i_{M'}^*TM').
\]
The difference of these equal-rank bundles has zero rational Pontryagin
character and hence belongs to $T$, by \cref{lem:KOCP} (and directly for $n=1$).
Equation~\eqref{eq:J-injective-torsion} shows that the difference is zero
in $K$, so the two top coefficients agree.  For a homeomorphism, the required
rational Pontryagin equality is supplied by topological invariance of rational
Pontryagin classes.  This proves the asserted homeomorphism invariance.
\end{proof}

\begin{proposition}\label{prop:delta-weights}
If \(n=4r+1\), then
\[
\delta(M_{a,b})=\Delta_r(a,b).
\]
\end{proposition}

\begin{proof}
Pulling the stable tangent bundle formula
\[
TM_{a,b}\oplus\varepsilon^3
\cong
f^*\Gamma_{a,b,\R}
\]
back along the map
$i:\CP^n\longrightarrow M_{a,b}$
gives
\begin{equation}\label{eq:pullback-tangent-delta}
i^*TM_{a,b}\oplus\varepsilon^3
\cong
j_n^*\Gamma_{a,b,\R}.
\end{equation}
Since adding a trivial bundle does not change a reduced \(KO\)-class,
the invariant \(\delta(M_{a,b})\) is the top Fujii coefficient of
\(j_n^*\Gamma_{a,b,\R}\).

We first compute this coefficient for a single summand
\((L^m)_\R\).

Put
\[
E_m=(L^m)_\R,\qquad A_m=[E_m]\in KO(\CP^n),
\]
where \(L\) is the canonical complex line bundle.  The real tensor
product of underlying real bundles satisfies
\[
(L)_\R\otimes_\R(L^m)_\R
\cong
(L^{m+1})_\R\oplus(L^{1-m})_\R.
\]
Indeed, for complex vector bundles \(E,F\),
\[
E_\R\otimes_\R F_\R
\cong
(E\otimes_\C F)_\R
\oplus
(E\otimes_\C\overline F)_\R.
\]
Taking \(E=L\) and \(F=L^m\) gives the displayed decomposition.
Moreover,
\[
L^{1-m}
=
\overline{L^{m-1}},
\]
and a complex bundle and its conjugate have isomorphic underlying real
bundles.  Hence
\[
(L^{1-m})_\R\cong(L^{m-1})_\R,
\]
so that
\[
E_1\otimes_\R E_m
\cong
E_{m+1}\oplus E_{m-1}.
\]
Passing to \(KO(\CP^n)\), we obtain the recurrence
\begin{equation}\label{eq:Am-recurrence}
A_{m+1}=A_1A_m-A_{m-1}.
\end{equation}

Now
$A_0=2,$ $A_1=[L_\R]=2+\mu$
where
$\mu=[L_\R]-2\in\widetilde{KO}(\CP^n).$
On the other hand, the Chebyshev polynomials of the first kind satisfy
\[
T_{m+1}(z)=2zT_m(z)-T_{m-1}(z),
\qquad
T_0(z)=1,\quad T_1(z)=z.
\]
The following expressions are polynomials in $\Z[t]$; they are
formed before evaluating at $\mu$.  In particular, no division by $2$ in
$KO(\CP^n)$ is intended.  The polynomials
\[
Q_m(t):=2T_m\left(1+\frac{t}{2}\right)
\]
satisfy
\[
Q_0(t)=2,\qquad Q_1(t)=2+t,
\]
and
\[
Q_{m+1}(t)=(2+t)Q_m(t)-Q_{m-1}(t).
\]
Therefore \(A_m\) and \(Q_m(\mu)\) have the same initial values and
satisfy the same recurrence.  By induction,
\[
[(L^m)_\R]=A_m=Q_m(\mu)
=
2T_m\left(1+\frac{\mu}{2}\right).
\]
Consequently
\begin{equation}\label{eq:real-line-Rm}
[(L^m)_\R]-2
=
R_m(\mu),
\end{equation}
where
\[
R_m(t)
:=
2T_m\left(1+\frac{t}{2}\right)-2.
\]

Write
\begin{equation}\label{eq:Rm-expansion}
R_m(t)
=
\sum_{s\ge1}c_s(m)t^s.
\end{equation}
The standard coefficient formula for the Chebyshev polynomial gives
\begin{equation}\label{eq:csm}
c_s(m)=
\begin{cases}
\displaystyle
\frac{m}{s}
\binom{m+s-1}{2s-1},
& m\ge s,\\[1.2ex]
0,
& m<s.
\end{cases}
\end{equation}
Substituting \(t=\mu\) in \eqref{eq:Rm-expansion}, we obtain
\begin{equation}\label{eq:real-line-expansion}
[(L^m)_\R]-2
=
\sum_{s\ge1}c_s(m)\mu^s \in \widetilde{KO}(\CP^n).
\end{equation}

We now use Fujii's description of the real \(K\)-group.  Since
\(n=4r+1\),
\[
\widetilde{KO}(\CP^n)
\cong
\Z\{\mu,\mu^2,\ldots,\mu^{2r}\}
\oplus
(\Z/2)\{\mu^{2r+1}\}.
\]
Moreover,
\[
\mu^{2r+2}=0.
\]
Hence every reduced \(KO\)-class has a unique expression
\[
q_1\mu+\cdots+q_{2r}\mu^{2r}
+
\varepsilon\mu^{2r+1},
\qquad
q_i\in\Z,\quad
\varepsilon\in\Z/2.
\]
By definition, \(\delta\) is precisely the last coefficient
\(\varepsilon\).

Applying this to \eqref{eq:real-line-expansion}, the coefficient of
\(\mu^{2r+1}\) is \(c_{2r+1}(m)\).  Since
\(\mu^{2r+1}\) has order two, only the parity of this integer matters.
Thus, for \(m\ge0\),
\begin{equation}\label{eq:delta-line-positive}
\delta((L^m)_\R)
=
c_{2r+1}(m)\pmod2.
\end{equation}

If \(m<0\), then
$L^{-m}\cong\overline{L^m},$
and \(L^m\) and \(\overline{L^m}\) have isomorphic underlying real
bundles.  Therefore
\[
(L^m)_\R\cong(L^{|m|})_\R,
\]
and hence, for every \(m\in\Z\),
\begin{equation}\label{eq:delta-line}
\delta((L^m)_\R)
=
c_{2r+1}(|m|)\pmod2.
\end{equation}

Finally,
\[
\Gamma_{a,b,\R}
=
\bigoplus_{i=1}^{n+1}(L^{a_i})_\R
\oplus
\bigoplus_{j=1}^{n+1}(L^{b_j})_\R,
\]
so
\[
\begin{aligned}
\delta(M_{a,b})
&=
\sum_{i=1}^{n+1}
c_{2r+1}(|a_i|)
+
\sum_{j=1}^{n+1}
c_{2r+1}(|b_j|)
\pmod2\\
&=
\Delta_r(a,b).
\end{aligned}
\]
This proves the claim.
\end{proof}

\section{The smooth normal type}\label{sec:normal-type}

Let \(M=M_{a,b}\), and let \(i:\CP^n\to M\) be the projective-stage map.
There is a cofibration
\begin{equation}\label{eq:cofiber}
\CP^n\xrightarrow{i}M\xrightarrow{q}Q_M,
\end{equation}
where \(Q_M\) is \(2n\)-connected and
\(H_{2n+1}(Q_M)\cong\Z\).

Let $c_r:\BO\langle r\rangle\to\BO$ denote the $(r-1)$-connected cover: its
source has zero homotopy groups below $r$, and $c_r$ induces an isomorphism in
degrees at least $r$.  We shall use $r=2n+2$.

From now on, stable normal classes and their extensions are taken to
have virtual rank zero: $\nu_M^0=[\eps^{4n+1}]-[TM]$.  We suppress the
superscript on $\nu_M$ when using reduced $KO$-theory.  This convention is
also used in the Thom spectra and their Atiyah--Hirzebruch filtrations.
Recall that $f\circ i\simeq j_n:\CP^n\to\CP^\infty$.  Let $\xi$ be a
virtual rank-zero bundle over $\CP^\infty$ such that
\begin{equation}\label{eq:xi-restriction}
 j_n^*\xi\cong_{st} i^*\nu_M,
\end{equation}
where $\nu_M$ is the stable normal bundle.  By exactness of reduced $KO$-theory
for \eqref{eq:cofiber}, there is a class
$\bar\rho_M\in\widetilde{KO}(Q_M)$ satisfying
\begin{equation}\label{eq:residual}
 q^*\bar\rho_M=\nu_M-f^*\xi.
\end{equation}

The next lemma is where the special geometry of the orbit manifold is used.

\begin{lemma}\label{lem:orbit-lift}
For every $n\ge2$, the classifying map of $\bar\rho_M$ lifts to
$\BO\langle2n+2\rangle$.
\end{lemma}

\begin{proof}
Let
\[
p:P=S^{2n+1}\times S^{2n+1}\longrightarrow M
\]
be the principal circle-bundle projection.  Since $P$ is stably parallelizable and
$TP\cong p^*TM\oplus\eps^1$, the pullback $p^*\nu_M$ is stably trivial.  The
composition $f\circ p$ is null-homotopic because $f$ classifies the principal
bundle.  Pulling back \eqref{eq:residual} therefore gives
\begin{equation}\label{eq:pullback-residual-zero}
 (q\circ p)^*\bar\rho_M=0\quad\text{in }\widetilde{KO}(P).
\end{equation}

We first show that
\[
(q\circ p)_*:H_{2n+1}(P)\longrightarrow H_{2n+1}(Q_M)
\]
is surjective.  The
Gysin sequence of $p$ shows that
$p^*:H^{2n+1}(M)\to H^{2n+1}(P)$ is a primitive injection.  Dually,
$p_*:H_{2n+1}(P)\to H_{2n+1}(M)$ is surjective, while
$q_*:H_{2n+1}(M)\to H_{2n+1}(Q_M)$ is an isomorphism.

The product \(P\) is \(2n\)-connected, so the Hurewicz map
$\pi_{2n+1}(P)\to H_{2n+1}(P)$ is an isomorphism.  Choose
$e:S^{2n+1}\to P$ so that $q\circ p\circ e$ represents a generator of
$H_{2n+1}(Q_M)$.  By \eqref{eq:pullback-residual-zero}, the restriction of
$\bar\rho_M$ to this bottom generator is stably trivial.

Since $Q_M$ is $2n$-connected, the only possible obstruction to lifting its
classifying map through the $(2n+1)$-connected cover
$\BO\langle2n+2\rangle\to\BO$ is the component detected on
$\pi_{2n+1}(Q_M)$.  We have just shown that this component is zero.  Hence the
lift exists.
\end{proof}

\begin{proposition}\label{prop:common-normal-type}
Let $M=M_{a,b}$ and $M'=M_{a',b'}$ be two orbit manifolds.  Suppose
\[
 \Pi_k(a,b)=\Pi_k(a',b')
 \quad(1\le k\le\lfloor n/2\rfloor),
\]
and, if $n=4r+1$, suppose also
$\Delta_r(a,b)=\Delta_r(a',b')$.  Then $M$ and $M'$ admit normal
$2n$-smoothings into a common normal $2n$-type
\begin{equation}\label{eq:normal-type}
 p_\xi:B_\xi=\CP^\infty\times\BO\langle2n+2\rangle\longrightarrow\BO.
\end{equation}
\end{proposition}

\begin{proof}
By \eqref{eq:pk}, equality of the $\Pi_k$ is equality of all Pontryagin classes
that can occur on the projective stage.  If $n\not\equiv1\pmod4$,
\cref{lem:KOCP} shows that $i^*\nu_M$ and $i'^*\nu_{M'}$ are stably isomorphic.
If $n=4r+1$, the same conclusion follows from \cref{prop:delta-weights} and
the additional equality of $\Delta_r$.  Passing from tangent to stable normal bundles does not change this
conclusion: the total Pontryagin class is inverted, while the top
\(\Z/2\)-coefficient changes sign and is therefore unchanged.  Choose a
virtual bundle \(\xi\) over \(\CP^\infty\) extending this common stable
class.

By \cref{lem:orbit-lift}, choose lifts
\[
 \mu_M:Q_M\to\BO\langle2n+2\rangle,
 \qquad
 \mu_{M'}:Q_{M'}\to\BO\langle2n+2\rangle
\]
of the residual normal bundles.  Define
\begin{equation}\label{canoNormalsmoothing}
	\bar\nu_M=(f,\mu_M\circ q):M\to B_\xi,
 \qquad
 \bar\nu_{M'}=(f',\mu_{M'}\circ q'):M'\to B_\xi.
\end{equation}
The map \(p_\xi\) classifies the Whitney sum of \(\xi\) with the
universal stable bundle on the connected-cover factor.  Hence the maps
in \eqref{canoNormalsmoothing} classify the stable normal bundles.
Because $f$ and $f'$ are $(2n+1)$-connected and the second factor of $B_\xi$
is $(2n+1)$-connected, both $\bar\nu_M$ and $\bar\nu_{M'}$ are
$(2n+1)$-equivalences.  They are therefore normal $2n$-smoothings into the
same normal $2n$-type.
\end{proof}

\section{The smooth normal bordism filtration}\label{sec:smooth-filtration}

Fix two orbit manifolds satisfying the hypotheses of
\cref{prop:common-normal-type}, and let $\xi$ be the common virtual bundle on
$\CP^\infty$.  The common normal $2n$-type gives bordism groups 
\[
 \Omega_*^{O\langle2n+2\rangle}(\xi)
 \cong \pi_*\bigl(M\xi\wedge MO\langle2n+2\rangle\bigr).
\]
Its homological Atiyah--Hirzebruch spectral sequence is
\begin{equation}\label{eq:AHSS-general}
 E^2_{s,t}=H_s\bigl(M\xi;\Omega_t^{O\langle2n+2\rangle}\bigr)
 \Longrightarrow \Omega_{s+t}^{O\langle2n+2\rangle}(\xi).
\end{equation}
For a fixed total degree $m$ let $F_{s,m-s}$ denote the skeletal filtration,
so that
\[
 F_{s,m-s}/F_{s-1,m-s+1}\cong E^\infty_{s,m-s}.
\]
We shall mainly use $m=4n+1$.

\begin{lemma}\label{lem:finite}
The group $\Omega_{4n+1}^{O\langle2n+2\rangle}(\xi)$ is finite.
\end{lemma}

\begin{proof}
The spectrum $MO\langle2n+2\rangle$ is connective and of finite type.
Rationally its homology is concentrated in degrees divisible by four, while
$H_*(M\xi;\mathbb Q)\cong H_*(\CP^\infty;\mathbb Q)$ is concentrated in even
degrees.  Hence every $E^2$-term in the odd total degree $4n+1$ is finite,
and so is the group $\Omega_{4n+1}^{O\langle2n+2\rangle}(\xi)$.
\end{proof}

\subsection{Equivariant prime avoidance}\label{sec:prime}

We now isolate the equivariant input in the bordism argument.  It comes
directly from the two coordinate projections of the linear action.

For a weight vector $c=(c_1,\ldots,c_{n+1})$, let $S_c^{2n+1}$ denote the
sphere $S^{2n+1}\subset\C^{n+1}$ with the weighted action
$z\cdot(w_1,\ldots,w_{n+1})=(z^{c_1}w_1,\ldots,z^{c_{n+1}}w_{n+1})$, and set
\[
 X_c=ES^1\times_{S^1}S_c^{2n+1}.
\]
The map $h_c:X_c\to\CP^\infty$ is the unit sphere bundle of
$\bigoplus_iL^{c_i}$.

\begin{lemma}\label{lem:Xc-homology}
Let $C(c)=\prod_i c_i$ and suppose $k=|C(c)|>0$.  Then
\[
 H^{2j}(X_c;\Z)\cong
 \begin{cases}
 \Z,&0\le j\le n,\\
 \Z/k,&j\ge n+1,
 \end{cases}
 \qquad
 H^{2j+1}(X_c;\Z)=0.
\]

\end{lemma}

\begin{proof}
The Euler class of the oriented sphere bundle $h_c$ is
$C(c)x^{n+1}$.  The Gysin sequence gives the displayed cohomology groups.
\end{proof}

\begin{proposition}[Prime-avoidance factorization]\label{prop:prime}
Let $M=M_{a,b}$ and let $f:M\to\CP^\infty$ be the map in
\eqref{eq:borel}.  For every prime $p$ there exist a space $X_p$ and maps
\[
 M\xrightarrow{g_p}X_p\xrightarrow{h_p}\CP^\infty
\]
such that $h_p\circ g_p\simeq f$ and
\[
 H_s(X_p;\Z_{(p)})=0\qquad(s\ge2n+1).
\]
\end{proposition}

\begin{proof}
The two coordinate projections
\[
 \operatorname{pr}_1:P\to S_a^{2n+1},
 \qquad
 \operatorname{pr}_2:P\to S_b^{2n+1}
\]
are $S^1$-equivariant for the actions defined.  They therefore induce
maps of Borel constructions
\[
 M\simeq ES^1\times_{S^1}P\longrightarrow X_a,
 \qquad
 M\longrightarrow X_b,
\]
commuting with the maps to $\CP^\infty$.

By \eqref{eq:gcdAB}, for every prime $p$ at least one of the nonzero integers
$|A(a)|$ and $|B(b)|$ is prime to $p$.  Choose the corresponding factor and
call its Borel construction $X_p$.  The Gysin calculation and the integral universal coefficient theorem give
\[
 H_s(X_p;\Z)=
 \begin{cases}
 \Z,&s=0,2,\ldots,2n,\\
 \Z/k,&s\ge2n+1\text{ odd},\\
 0,&\text{otherwise},
 \end{cases}
\]
where $k$ is the absolute value of the chosen product of weights.
Since $p\nmid k$, localizing at $p$ proves the assertion.
\end{proof}

\begin{remark}\label{rem:scope}
The proof of \cref{prop:prime} uses the actual equivariant coordinate
projection of the linear action.  No equivariant projection is needed or
asserted for an arbitrary principal circle bundle whose total space is merely
homotopy equivalent to a product of spheres.  This distinction is essential
to the scope of the present theorem.
\end{remark}

\begin{proposition}\label{prop:low-filtration}
For every \(n\ge2\), the canonical normal smoothing
\eqref{canoNormalsmoothing} of \(M=M_{a,b}\) satisfies
\[
 [M,\bar\nu_M]\in F_{2n,2n+1}
 \subset \Omega_{4n+1}^{O\langle2n+2\rangle}(\xi).
\]
\end{proposition}

\begin{proof}
Fix a prime $p$ and the factorization
$M\xrightarrow{g_p}X_p\xrightarrow{h_p}\CP^\infty$ of
\cref{prop:prime}.  Put $\xi_p=h_p^*\xi$.
The homotopy $h_p g_p\simeq f$ transports the stable normal bundle
identification, while the same connected-cover component supplies a
normal bordism class over $X_p$ mapping to $[M,\bar\nu_M]$.

For every finitely generated $\Z_{(p)}$-module $A$, the homological
universal coefficient sequence and \cref{prop:prime} give
\[
 H_s(X_p;A)=0\qquad(s\ge2n+1).
\]
At $s=2n+1$, the possible Tor term also vanishes because
$H_{2n}(X_p;\Z_{(p)})\cong\Z_{(p)}$ is free.  The virtual bundle
$\xi_p$ is oriented and has rank zero, so the same assertion holds for
the homology of its Thom spectrum.  Consequently the localized AHSS
for bordism over $X_p$ has no terms above filtration $2n$ in the
relevant total degree.  Naturality puts $[M,\bar\nu_M]_{(p)}$ in that
filtration.  This holds for every prime.  Since the target bordism group
is finite by \cref{lem:finite}, the integral class belongs to
$F_{2n,2n+1}$ as claimed.
\end{proof}

\section{Passage to the PL normal type}\label{sec:PL}

Let \(\BPL\langle2n+2\rangle\) be the \((2n+1)\)-connected cover of
\(\BPL\), and let \(\MPL\langle2n+2\rangle\) be its Thom spectrum.
Forgetting the vector-bundle structure induces a map of normal types
\[
\id\times T:
\CP^\infty\times\BO\langle2n+2\rangle
\longrightarrow
\CP^\infty\times\BPL\langle2n+2\rangle
\]
and hence a homomorphism
\[
T_*:
\Omega_*^{O\langle2n+2\rangle}(\xi)
\longrightarrow
\Omega_*^{\PL\langle2n+2\rangle}(|\xi|).
\]
The target has an Atiyah--Hirzebruch spectral sequence with the same
base Thom spectrum \(M\xi\) and coefficient groups
\(\pi_*(\MPL\langle2n+2\rangle)\).

The following vanishing statement is the only point at which the
low-dimensional cases require separate attention.

\begin{lemma}\label{lem:PLcoeff}
For \(n\ge4\) and \(n\le i\le2n\),
\[
\pi_{2i+1}\bigl(\MPL\langle2n+2\rangle\bigr)=0.
\]
In addition,
\[
\pi_{11}(\MPL\langle8\rangle)
=
\pi_{13}(\MPL\langle8\rangle)
=0.
\]
\end{lemma}

\begin{proof}
For \(n\ge4\), this is the standard surgery calculation of highly
connected PL bordism.  A representative of
\(\Omega_{2i+1}^{\PL\langle2n+2\rangle}\) may be surgered below the
middle dimension, compatibly with its normal structure, until it is
\(i\)-connected.  It is then a PL homotopy sphere.  In dimensions at
least five, a PL homotopy sphere is PL homeomorphic to the standard
sphere, and its stable PL normal bundle is trivial.  If \(F\) denotes
the homotopy fiber of
\[
\BPL\langle2n+2\rangle\longrightarrow\BPL,
\]
then \(\pi_m(F)=0\) for \(m\ge2n+1\).  Hence the lifted
\(\BPL\langle2n+2\rangle\)-structure is null-homotopic, so the class
bounds.  This is the usual argument based on the first stability
theorem for PL bundles; see \cite[p.~247]{KirbySiebenmann1977}.
The same argument in the two low degrees gives the final displayed
vanishing statements.
\end{proof}

\begin{proposition}\label{MPL9}
The forgetful homomorphism
$
T_\#:\Omega_9^{O\langle8\rangle}
\longrightarrow
\Omega_9^{\PL\langle8\rangle}$
is zero.
\end{proposition}

\begin{proof}
By \cite{KervaireMilnor1963}, every class in
\(\Omega_9^{O\langle8\rangle}\) may be represented by a homotopy
\(9\)-sphere with its trivial \(\BO\langle8\rangle\)-structure.
After forgetting the smooth structure, such a sphere is PL
homeomorphic to \(S^9\), and the induced
\(\BPL\langle8\rangle\)-structure is trivial.  Its image therefore
vanishes in \(\Omega_9^{\PL\langle8\rangle}\).
\end{proof}

\begin{proposition}\label{prop:PLfiltration-zero}
For \(n\ge2\), the restriction $T_*|_{F_{2n,2n+1}}:
F_{2n,2n+1}
\longrightarrow
\Omega_{4n+1}^{\PL\langle2n+2\rangle}(|\xi|)$
is zero.
\end{proposition}

\begin{proof}
For \(s\le2n\), the homology of \(M\xi\) is concentrated in even
filtration degrees.  Thus a potentially nonzero term of total degree
\(4n+1\) has odd coefficient degree
\[
t=4n+1-s=2i+1,
\qquad n\le i\le2n.
\]
If \(n\ge4\), all such coefficient groups vanish by
\cref{lem:PLcoeff}.  Hence
$
F^{\PL}_{2n,2n+1}=0,$
and the desired result follows.

Suppose \(n=3\).  Since
\(\pi_7(\MO\langle8\rangle)=0\), one has
\(E^2_{6,7}=0\), and therefore
\begin{equation}\label{F67F49_13}
	F_{6,7}=F_{4,9}.
\end{equation}
On the PL side, \cref{lem:PLcoeff} gives
\(F^{\PL}_{3,10}=0\), so \(F^{\PL}_{4,9}\) is identified with the
\(E^\infty_{4,9}\)-term.  By \cref{MPL9}, the coefficient map in
degree \(9\) is zero; hence the induced maps on both
\(E^2_{4,9}\) and \(E^\infty_{4,9}\) are zero.  It follows that
\(T_*\) vanishes on \(F_{4,9}=F_{6,7}\).

Finally, let \(n=2\).  Since
$\BO\langle6\rangle\simeq\BO\langle8\rangle$, $\MO\langle8\rangle\simeq \MO\langle 6\rangle$. The relevant
coefficient groups are 
$
\pi_5(\MO\langle8\rangle)=\pi_7(\MO\langle8\rangle)=0;$
see \cite{Giambalvo1971}.  It follows that
$F_{4,5}=F_{0,9}.$
So \cref{MPL9} proves the assertion.
\end{proof}

The smooth normal \(2n\)-smoothing \eqref{canoNormalsmoothing} induces
the canonical PL normal \(2n\)-smoothing
\[
(f,T\circ\mu_M\circ q):
M_{a,b}
\longrightarrow
\CP^\infty\times\BPL\langle2n+2\rangle.
\]

\begin{corollary}\label{prop:PLbordant}
Under the hypotheses of \cref{prop:common-normal-type}, the canonical
PL normal \(2n\)-smoothings of \(M_{a,b}\) and \(M_{a',b'}\) are
normally bordant.
\end{corollary}

\begin{proof}
By \cref{prop:low-filtration}, both smooth classes lie in
\(F_{2n,2n+1}\).  Naturality carries them to
\(F^{\PL}_{2n,2n+1}\), and \cref{prop:PLfiltration-zero} shows that
both images vanish.
\end{proof}

\section{Modified surgery and the proof of the main theorem}\label{sec:surgery}

We now pass from PL normal bordism to PL homeomorphism.  We use
Kreck's original formulation of modified surgery
\cite[Theorem~5.2(b)]{KreckEarly}, with the extension to the PL
category described at the beginning of Section~1 of that manuscript;
compare \cite{Kreck1999}.  In the present simply connected situation,
the obstruction is controlled by a quotient of an ordinary Wall group.

For \(\mathrm{CAT}=O\) or \(\PL\), write
\[
 p_\xi^{\mathrm{CAT}}:B_\xi^{\mathrm{CAT}}
 =\CP^\infty\times \mathrm{BCAT}\langle2n+2\rangle
 \longrightarrow \mathrm{BCAT}
\]
for the map classifying the Whitney sum of the pullback of \(\xi\)
and the universal stable bundle on the connected-cover factor.
We replace this map by a fibration when necessary.  Here \(\xi\)
is regarded as a stable PL bundle when \(\mathrm{CAT}=\PL\), and \(\mathrm{CAT}=O\)
denotes the smooth category.

\begin{proposition}\label{prop:kreck-surgery} \quad

\begin{enumerate}[label=\textup{(\roman*)}]
\item
Let \(n\ge2\).  If two normal \(2n\)-smoothings
\[
 (N_0,\bar\nu_0),\ (N_1,\bar\nu_1)
 \longrightarrow B_\xi^{\PL}
\]
of closed simply connected \((4n+1)\)-manifolds are normally bordant,
then there exists a $\PL$ \(h\)-cobordism between \(N_0\) and \(N_1\)
carrying a \(B_\xi^{\PL}\)-normal structure extending the prescribed
boundary normal structures.

\item
When \(n=3\), let \(\mathrm{CAT}=O\) or \(\PL\).  Every
\(B_\xi^{\mathrm{CAT}}\)-normal bordism
\[
 (W^{14};N_0,N_1)\longrightarrow B_\xi^{\mathrm{CAT}}
\]
between normal \(6\)-smoothings of closed simply connected
\(13\)-manifolds is normally bordant relative to its boundary to
a \(\mathrm{CAT}\) \(h\)-cobordism.
\end{enumerate}
\end{proposition}

\begin{proof}
Put \(m=2n+1\), so that the bordism dimension is \(2m\) and the
boundary maps are normal \((m-1)\)-smoothings.  Thus
\cite[Theorem~5.2(b)]{KreckEarly} applies with \(k=m\).
The spaces \(B_\xi^{\mathrm{CAT}}\) have CW models of finite type, so the
finite-\((m+1)\)-skeleton hypothesis is satisfied.  As in that
theorem, we first perform compatible surgeries below the middle
dimension so that the normal map of the bordism is an
\(m\)-equivalence, without changing its relative normal bordism
class or its boundary data.  We also have
\[
 \pi_1(B_\xi^{\mathrm{CAT}})=1,
 \qquad w_1(B_\xi^{\mathrm{CAT}})=0,
 \qquad \operatorname{Wh}(1)=0.
\]

For a normal bordism \(V\) as above, Kreck's theorem gives a
relative normal bordism invariant \(\Theta(V)\) and the criterion
\begin{equation}\label{eq:kreck-criterion}
 \Theta(V)=0
 \quad\Longleftrightarrow\quad
 \begin{gathered}
 V\text{ is normally bordant relative to }\partial V\\
 \text{to a }\mathrm{CAT}\text{ }h\text{-cobordism}.
 \end{gathered}
\end{equation}
For \(m\ne3,7\), the obstruction takes values in
\[
 L_{2m}^{s,\tau}(\pi_1,w_1,S(V))=L_{2m}^{s,\tau}(1,0,S(V)).
\]
Since \(k=m\), the form parameter depends only on the boundary
normal smoothings:
\[
 S(V)=S(N_0\times I)=S(N_1\times I);
\]
see \cite[pp.~54--56]{KreckEarly}.  In particular, replacing the
bordism while preserving these boundary data does not change the
obstruction group.

\smallskip
\noindent\emph{The exceptional case \(n=3\).}
Here \(m=7\).  We first prove part~\textup{(ii)}, in either category.
The exceptional clause following
\cite[Theorem~5.2(b), p.~57]{KreckEarly} replaces the form parameter
by \(S(V)+\Z\) when
\[
 \bigl\langle
 (p_\xi^{\mathrm{CAT}})^*w_8,\pi_8(B_\xi^{\mathrm{CAT}})
 \bigr\rangle\ne0.
\]
Because \(\pi_8(\CP^\infty)=0\) and
\(\mathrm{BCAT}\langle8\rangle\to \mathrm{BCAT}\) induces an isomorphism on
\(\pi_8\), the homomorphism
\[
 (p_\xi^{\mathrm{CAT}})_*:
 \pi_8(B_\xi^{\mathrm{CAT}})\longrightarrow\pi_8(\mathrm{BCAT})
\]
is an isomorphism.  The octonionic Hopf \(8\)-plane bundle over
\(S^8\) has Euler number \(1\), hence nonzero \(w_8\).  Its class,
or its underlying PL bundle, therefore gives the required nonzero
pairing.  This is also the exceptional framing phenomenon in
\cite[Proposition~6]{Kreck1999}.

The group ring is \(\Z\), so the enlarged form parameter is
\[
 S(V)+\Z=\Z.
\]
The quadratic refinement consequently has zero target.  The remaining
forms are unimodular alternating forms over \(\Z\), all of which
are hyperbolic.  Since there is no Whitehead-torsion contribution,
\[
 L_{14}^{s,\tau}(1,0,\Z)=0;
\]
compare \cite[Example~4.10]{KreckEarly}.
Thus \(\Theta(V)=0\), and \eqref{eq:kreck-criterion} proves
part~\textup{(ii)}.  Taking \(\mathrm{CAT}=\PL\) also proves
part~\textup{(i)} when \(n=3\).

\smallskip
\noindent\emph{The $\PL$ case \(n\ne3\).}
Choose a \(B_\xi^{\PL}\)-normal bordism \(W\) between the prescribed
normal smoothings, and write \(S=S(W)\).  By
\cite[Lemma~4.1 and Section~4]{KreckEarly}, there is a natural
surjection
\[
 q_S:L_{4n+2}^{s}(\Z)
 \cong L_{2m}^{s,\tau}(1,0)
 \longrightarrow L_{2m}^{s,\tau}(1,0,S).
\]
Here we used \(\operatorname{Wh}(1)=0\).  Moreover,
\[
 L_{4n+2}^{s}(\Z)\cong\Z/2;
\]
see \cite[Example~4.10]{KreckEarly}.  Choose \(\alpha\) such that
\[
 q_S(\alpha)=-\Theta(W).
\]
If \(\Theta(W)=0\), choose \(\alpha=0\) and set \(W_\alpha=W\).

Suppose that \(\Theta(W)\ne0\), so \(\alpha\) is the nonzero Wall
class.  We realize this correction by a PL construction, rather
than by a normal bordism of the given \(W\).  Let \(P^{2m}\) be the
parallelizable Kervaire--Arf plumbing of two tangent disk bundles
of \(S^m\).  Its boundary is a homotopy \((2m-1)\)-sphere and hence
is PL homeomorphic to the standard sphere.  Capping this boundary
by a PL disk gives a closed PL manifold \(K^{2m}\); see
\cite[Section~1]{CrowleyHambleton2015}.

The stable PL normal bundle of \(K\) is trivial on \(P\).
Choose the trivialization supplied by the plumbing.  Relative to this
trivialization, its classifying map factors through
\[
 K/P\cong S^{2m}.
\]
Since \(\BPL\langle m+1\rangle\to \BPL\) induces an isomorphism on
\(\pi_{2m}\), this factorization lifts to
\(\BPL\langle m+1\rangle\), compatibly with the chosen
trivialization on \(P\).  Together with the constant map to
\(\CP^\infty\), the lift gives a \(B_\xi^{\PL}\)-normal
structure on \(K\).  Notice that this only uses triviality of the
stable normal bundle on the plumbing; no global stable
trivialization of the closed PL manifold \(K\) is asserted.

Form the interior connected sum
\[
 W_\alpha=W\#_{\mathrm{int}}K.
\]
Choose the summing disks and their normal trivializations compatibly.
The normal structures then glue, and the prescribed normal structures
on \(N_0\) and \(N_1\) are unchanged.

We check the effect on the obstruction.  The two plumbing cores give
classes \(e,f\in H_m(K;\Z)\) with
\[
 \lambda_K(e,f)=1,
 \qquad \mu_K(e)=\mu_K(f)=1\pmod2.
\]
Thus their quadratic form has Arf invariant one and represents
\(\alpha\).  The normal map of \(K\) is constant on the plumbing
spine, so these classes lie in the surgery kernel.  Under interior
connected sum, representatives from \(W\) and the plumbing may be
chosen disjoint.  Their mixed intersections vanish, and the boundary
images lie entirely in the \(W\)-summand.  Consequently the reduced
middle-dimensional form is the orthogonal sum of the original form
and the plumbing form, with quadratic values reduced to \(\Z/S\).
In Kreck's obstruction group this gives
\[
 \Theta(W_\alpha)=\Theta(W)+q_S(\alpha)=0.
\]
The form parameter remains \(S\), since the boundary normal
smoothings have not changed.  Applying \eqref{eq:kreck-criterion}
to \(W_\alpha\) now yields the required PL \(h\)-cobordism,
proving part~\textup{(i)}.
\end{proof}

\begin{proof}[Proof of \cref{thm:main}]
Suppose first that \(n=1\).  The cohomology computation and
\cref{rem:n1-w2} give
\[
 H_2(M;\Z)\cong\Z,
 \qquad
 w_2(M)=\Delta_0\bar x.
\]
The Smale--Barden classification \cite{Smale1962,Barden1965} therefore gives exactly two
diffeomorphism types, distinguished by \(\Delta_0\).

Assume now \(n\ge2\).  Necessity follows from the characteristic-class
calculation.  A homeomorphism sends the degree-two generator to its
positive or negative.  Novikov invariance of the rational Pontryagin
classes, together with torsion-freeness of the relevant integral
cohomology groups, gives
\[
 \Pi_k(a,b)=\Pi_k(a',b')
 \qquad
 (1\le k\le\lfloor n/2\rfloor).
\]
When \(n=4r+1\), invariance of the top \(KO\)-coefficient follows from
\cref{prop:delta-intrinsic}, and \cref{prop:delta-weights} gives
\[
 \Delta_r(a,b)=\Delta_r(a',b').
\]

Conversely, assume that these arithmetic invariants agree.
By \cref{prop:common-normal-type}, the two manifolds admit normal
\(2n\)-smoothings into a common smooth normal type, and by
\cref{prop:PLbordant} the induced PL normal smoothings are normally
bordant.  Applying \cref{prop:kreck-surgery} with \(\mathrm{CAT}=\PL\) changes
this bordism, relative to its boundary, into a PL \(h\)-cobordism.
Since the manifolds are simply connected, the PL \(h\)-cobordism
theorem gives
$ M_{a,b}\cong_{\PL}M_{a',b'}.$
\end{proof}

\section{The smooth classification in dimension thirteen}\label{sec:13smooth}

Throughout this section,
$
 M=M_{a,b},$
 $\lambda=\sum_{i=1}^{4}(a_i^2+b_i^2).$
By \cref{prop:tangent},
\begin{equation}\label{eq:13-p1w2}
 p_1(M)=\lambda x^2,
 \qquad
 w_2(M)\equiv\lambda x\pmod2.
\end{equation}

\subsection{Non-spin case}\label{subsec:13-nonspin}
Fix a value $\lambda$ and let $\xi_\lambda$ be the stable bundle on $\CP^\infty$ in the common normal $6$-type.  
We choose this extension
canonically as follows.  If
\[
 u=[L_{\mathbb R}]-2\in\widetilde{KO}(\CP^\infty),
\]
then we take $\xi_\lambda=-\lambda u$.  Its restriction to $\CP^3$ is the
required stable normal class, since $p_1$ detects
$\widetilde{KO}(\CP^3)$.  In particular,
\begin{equation}\label{eq:13-xi-char}
 p_1(\xi_\lambda)=-\lambda x^2,
 \qquad
 w_2(\xi_\lambda)\equiv\lambda x\pmod2.
\end{equation}

The general filtration estimate can be sharpened slightly in this
dimension.

\begin{lemma}\label{lem:13-F4}
For every \(13\)-dimensional orbit manifold \(M\), the canonical normal
\(6\)-smoothing satisfies
$[M,\bar\nu_M]\in F_{4,9}\subset \Omega_{13}^{O\langle 8\rangle}(\xi_\lambda).$
\end{lemma}

\begin{proof}
This lemma has been proved (see \eqref{F67F49_13}).
\end{proof}

We shall use the classical low-dimensional groups
\begin{equation}\label{eq:string-low}
 \Omega_8^{O\langle8\rangle}\cong\mathbb Z\oplus\mathbb Z/2,
 \quad
 \Omega_9^{O\langle8\rangle}\cong(\mathbb Z/2)^2,
 \quad
 \Omega_{10}^{O\langle8\rangle}\cong\mathbb Z/6,
 \quad \Omega_{11}^{O\langle8\rangle}\cong 0,
 \quad
 \Omega_{13}^{O\langle8\rangle}\cong\mathbb Z/3,
\end{equation}
see \cite{Giambalvo1971,HoveyRavenel1995}.

\begin{lemma}\label{lem:13-differentials}
In the $\mathrm{AHSS}$ for
$\Omega_*^{O\langle8\rangle}(\xi_\lambda)$,
$d_2:E^2_{6,8}\to E^2_{4,9}$ is mod-$2$ reduction followed by the
homological dual of
\[
 Sq^2:H^4(M\xi_\lambda;\mathbb Z/2)
 \longrightarrow H^6(M\xi_\lambda;\mathbb Z/2).
\]
\end{lemma}

\begin{proof}
The basic $d_2$ is the low-dimensional homological AHSS
formula of \cite{Shen2026Differentials}.  Recall that
$ \Omega_8^{O\langle8\rangle}
 =\mathbb Z\langle\omega\rangle\oplus
  \mathbb Z/2\langle c_0\rangle$
so that
$
 \Omega_9^{O\langle8\rangle}
 =\mathbb Z/2\langle\eta\omega\rangle\oplus
  \mathbb Z/2\langle\eta c_0\rangle.$
Multiplicativity identifies the two components of
$d_2:E^2_{6,8}\to E^2_{4,9}$ with the basic $d_2$ tensored respectively with
$\omega$ and $c_0$.  This yields the statement.
\end{proof}

%
%
%

\begin{proposition}\label{prop:13-F4-compute}
Assume that $\lambda$ is odd.  Then $F_{4,9}$ is isomorphic to $0$ or $\Z/3$. If $F_{4,9}\cong \Z/3$, its generator is represented by an exotic $13$-sphere with its canonical
$O\langle8\rangle$ normal structure.
\end{proposition}

\begin{proof}
Since $\lambda$ is odd, $w_2(\xi_\lambda)\ne0$.  If $U$ is the Thom class,
then
\[
 Sq^2(U)=w_2(\xi_\lambda)U.
\]
Because $H^*(M\xi_\lambda;\mathbb Z/2)$ is a free rank-one module over
$H^*(\CP^\infty;\mathbb Z/2)$ on $U$, the operation
\[
 Sq^2:H^4(M\xi_\lambda;\mathbb Z/2)
 \longrightarrow H^6(M\xi_\lambda;\mathbb Z/2)
\]
is an isomorphism.  By \cref{lem:13-differentials}, the differential
$d_2:E^2_{6,8}\to E^2_{4,9}$ is onto and kills the whole
$E^2_{4,9}\cong(\mathbb Z/2)^2$.  The terms in filtrations $1,2,3$ in total
degree $13$ vanish, so
\[
 F_{4,9}=F_{0,13}=E^\infty_{0,13}.
\]

Since
$E^2_{0,13}=\Omega_{13}^{O\langle8\rangle}\cong\mathbb Z/3$, $E^\infty_{0,13}$ is $0$ or $\Z/3$. If $E^\infty_{0,13}\cong \Z/3$, the edge homomorphism of the AHSS gives the isomorphism 
$$\Omega_{13}^{O\langle8\rangle}\cong E^\infty_{0,13}.$$
By \cite{KervaireMilnor1963}, the nonzero element of $\Omega_{13}^{O\langle8\rangle}$ may be represented by an
exotic $13$-sphere; in the present dimension its stable normal bundle is
trivial, and the resulting $O\langle8\rangle$ normal structure is the
canonical one. This finishes the proof.
\end{proof}

\begin{theorem}\label{thm:13-smooth}
Let $M=M_{a,b}$ and $M'=M_{a',b'}$ be $13$-dimensional orbit manifolds with
$\lambda=\lambda'\equiv 1\pmod 2$.
 Then, there is a homotopy sphere
$\Sigma^{13}\in\Theta_{13}\cong\mathbb Z/3$ such that
\[
 M\cong_{\mathrm{Diff}}M'\#\Sigma^{13}.
\]
\end{theorem}

\begin{proof}
The two manifolds have the same smooth normal $6$-type by
\cref{prop:common-normal-type}.  By \cref{lem:13-F4}, their normal smoothing
classes lie in $F_{4,9}$.  If $\lambda$ is odd,
\cref{prop:13-F4-compute} shows that their difference is represented by a homotopy $13$-sphere. Thus for some $\Sigma^{13}$ the normal
smoothings of $M$ and $M'\#\Sigma^{13}$ are bordant.  Apply
\cref{prop:kreck-surgery} in the smooth category to obtain the diffeomorphism.
\end{proof}

Next we remove the homotopy sphere for linear quotients.

\begin{lemma}\label{lem:orientation}
Every orbit manifold $M_{a,b}$ in this section has a self-diffeomorphism
that reverses orientation.
\end{lemma}

\begin{proof}
Simultaneous complex conjugation
\[
 R:S^7\times S^7\longrightarrow S^7\times S^7,
 \qquad
 R(u,v)=(\bar u,\bar v),
\]
satisfies
\[
 R(t\cdot(u,v))=t^{-1}\cdot R(u,v).
\]
It therefore descends to a self-diffeomorphism $r_M$ of the quotient.  Complex
conjugation on $S^7\subset\mathbb C^4$ has degree $(-1)^4=+1$, so $R$ preserves the
orientation of $S^7\times S^7$.  On the oriented circle fibre it induces inversion
$t\mapsto t^{-1}$, which has degree $-1$.  With the quotient orientation convention,
$r_M$ therefore has degree $-1$.
\end{proof}

\begin{lemma}\label{lem:odd-sphere-cancel}
Let $X$ and $Y$ be closed oriented manifolds, each admitting an
orientation-reversing self-diffeomorphism.  Suppose
\[
 X\cong_+Y\#\Sigma
\]
for a homotopy sphere $\Sigma$ whose class has odd order in the group of homotopy
spheres.  Then $X\cong_+Y$.
\end{lemma}

\begin{proof}
Reverse orientations in the displayed diffeomorphism and then use the
orientation-reversing self-diffeomorphisms of $X$ and $Y$.  This gives
\[
 X\cong_+Y\#(-\Sigma).
\]
Comparing the two expressions yields
\[
 Y\#\Sigma\cong_+Y\#(-\Sigma).
\]
Taking connected sum with $\Sigma$ gives
\[
 Y\#(2\Sigma)\cong_+Y,
\]
so \(2[\Sigma]\) belongs to the oriented inertia group \(I(Y)\).
The cyclic subgroup generated by \([\Sigma]\) has odd order, and hence
multiplication by \(2\) is an automorphism on that subgroup.  Therefore
\([\Sigma]\in I(Y)\), and consequently \(Y\#\Sigma\cong_+Y\).
\end{proof}

\begin{theorem}[Non-spin smooth classification in dimension \(13\)]
\label{thm:13-complete-diff}
Let \(M_{a,b}\) and \(M_{a',b'}\) be non-spin \(13\)-dimensional orbit
manifolds.  Then
\[
M_{a,b}\cong_{\Diff}M_{a',b'}
\quad\Longleftrightarrow\quad
\lambda(a,b)=\lambda(a',b').
\]
Moreover, the diffeomorphism may be chosen orientation-preserving.
\end{theorem}

\begin{proof}
Necessity follows from
\(p_1(M_{a,b})=\lambda(a,b)x^2\).
Conversely, suppose that the two \(\lambda\)-invariants agree.  Since
the manifolds are non-spin, their common value is odd.
By \cref{thm:13-smooth}, there is a homotopy sphere
\(\Sigma^{13}\in\Theta_{13}\cong\mathbb Z/3\) such that
\[
M_{a,b}\cong_+M_{a',b'}\#\Sigma^{13}.
\]
Both orbit manifolds admit orientation-reversing self-diffeomorphisms
by \cref{lem:orientation}.  Since every element of
\(\Theta_{13}\cong\mathbb Z/3\) has odd order,
\cref{lem:odd-sphere-cancel} removes the homotopy sphere.
\end{proof}

\subsection{\texorpdfstring{The case of $\lambda\equiv0\pmod4$}{The case lambda congruent to 0 modulo 4}}
\label{sec:spin4}

Throughout this subsection, we use generators
\[
 H^*(M;\mathbb Z)=\mathbb Z[x,y]/(x^4,y^2),
 \qquad |x|=2,\quad |y|=7,
 \qquad \langle x^3y,[M]\rangle=1.
\]
A bar denotes reduction modulo two.  We choose a finite CW model with
a single top cell and use the corresponding skeletal filtration.
We write
\[
\mathcal F^s[M,G/O]
=
\ker\bigl([M,G/O]\longrightarrow[M^{(s-1)},G/O]\bigr)
\]
for the skeletal filtration associated to the chosen CW model.

The aim of this subsection is to prove that 
\begin{theorem}\label{thm:spin4}
If $\lambda(a,b)=\lambda(a',b')\equiv 0\pmod 4$, then $M_{a,b}\cong_{\Diff}M_{a',b'}$
\end{theorem}

  The proof starts with the PL classification:
orbit manifolds with the same first Pontryagin coefficient are
orientation-preservingly PL homeomorphic.  We then show that the
smooth normal invariant of such a PL homeomorphism can be removed.

Let $M=M_{a,b}$, $M'=M_{a',b'}$, $\lambda=\lambda(a,b)$. We first show
\begin{lemma}\label{Normal_invariant_forh}
	Let $h:M'\to M$ be an orientation-preservingly PL homeomorphic. Then its smooth normal invariant satisfies $\eta^{\Diff}(h)\in \mathcal F^9[M,G/O].$
\end{lemma}

The main point is to eliminate the possible \(2\)-primary obstruction of $h$
in filtration nine.  We construct a tangential self-homotopy
equivalence:
\begin{lemma}\label{New_self_homotopy}
	There exists a self-homotopy
equivalence \(f_{\nu^2}^M\) whose smooth normal invariant is
the nonzero element
$
\theta_M\in\mathcal F^9[M,G/O]_{(2)}\cong\mathbb Z/2.$
\end{lemma}

 After the \(2\)-primary component of the smooth normal invariant has been eliminated, the remaining normal invariant lies entirely in the top filtration and can be absorbed by a \(13\)-dimensional homotopy sphere.
 
 \begin{lemma}\label{3-primary_normal_in}
 	If there exists a homotopy
equivalence $g:M'\to M$ such that $\eta^{\Diff}(g)\in \mathcal F^9[M,G/O]$ and $\eta^{\Diff}(g)_{(2)}=0$, then $M\cong_{\Diff} M'\# \Sigma^{13}$ for some homotopy sphere $\Sigma^{13}$.
 \end{lemma}

We now prove the diffeomorphism classification result: 


\begin{proof}[Proof of Theorem \ref{thm:spin4}]
  The PL classification gives
a PL homeomorphism \(M'\to M\).  By \cref{lem:orientation}, we may
choose it orientation-preserving; denote it by \(h\).  By Lemma \ref{Normal_invariant_forh}, its smooth normal invariant belongs to
\(\mathcal F^9[M,G/O]\).

If $\eta^{\Diff}(h)=0\in \mathcal F^9[M,G/O]_{(2)}$, then $M\cong_{\Diff} M'\# \Sigma^{13}$ by Lemma \ref{3-primary_normal_in}.

If $0\ne \eta^{\Diff}(h)\in \mathcal F^9[M,G/O]_{(2)}$, we consider \(f_{\nu^2}^M\circ h\).  

Put \(f=f_{\nu^2}^M\). By cellular approximation, pullback by \(f^{-1}\) preserves the
skeletal filtration
\[
(f^{-1})^*
\bigl(\mathcal F^9[M,G/O]\bigr)
\subset
\mathcal F^9[M,G/O].
\]
On the \(2\)-primary subgroup
\(\mathcal F^9[M,G/O]_{(2)}\cong\mathbb Z/2\), this pullback is an
automorphism and therefore fixes the unique nonzero element
\(\theta_M\).
By Lemma \ref{New_self_homotopy} and the formula 
\[
\eta^{\mathrm{Diff}}(f\circ h)
=
\eta^{\mathrm{Diff}}(f)
+
(f^{-1})^*\eta^{\mathrm{Diff}}(h);
\]
see \cite[Lemma~2.5]{MadsenTaylorWilliams1980}, we have 
$$\eta^{\mathrm{Diff}}(f\circ h)\in \mathcal F^9[M,G/O]\quad\text{and} \quad
\eta^{\mathrm{Diff}}(f\circ h)_{(2)}=0.$$
By Lemma \ref{3-primary_normal_in}, $M\cong_{\Diff} M'\# \Sigma^{13}$ .
 
 Finally,
\cref{lem:odd-sphere-cancel} therefore yields
$M'\cong_{\Diff}M.$
\end{proof}

We now establish the ingredients used in the proof above.  We begin
with the Atiyah--Hirzebruch spectral sequence. In this subsection, we use the Whitney-sum group structure on
$[M,G/O]$.  The cohomological Atiyah--Hirzebruch spectral sequence is
\[
 E_2^{p,q}=H^p\bigl(M;\pi_{-q}(G/O)\bigr),
 \qquad
 d_r:E_r^{p,q}\longrightarrow E_r^{p+r,q-r+1}.
\]
For $p+q=0$, we have 
$$E_\infty^{p,q}\cong \mathcal F^p[M,G/O]/\mathcal F^{p+1}[M,G/O].$$

%

\subsubsection{The filtration-nine part of the normal invariant} $\quad$

Let \(j:SG\to G/O\) be the canonical map, where
\(SG=(QS^0)_1\), and let \(\nu\) denote a generator of
\(\pi_3^S\cong\pi_3(SG)\).  Write
\[
 a_9=j_*(\nu^3)=[\nu^3],\qquad b_9=j_*(\mu_9),
\]
where $\mu_9\in\pi_9^S$ is the classical stable homotopy element;
it is not in the image of the stable $J$-homomorphism.  The required
low-dimensional computations are
\begin{align}
 \pi_1(G/O)=\pi_3(G/O)=\pi_5(G/O)=\pi_7(G/O)&=0,\notag\\
 \pi_9(G/O)_{(2)}
 &=\mathbb Z/2\{a_9\}\oplus\mathbb Z/2\{b_9\},
 \label{eq:GO9}\\
 \pi_{10}(G/O)_{(2)}
 &=\mathbb Z/2\{j_*(\eta_\ast\mu_9)\},\notag\\
 \pi_{11}(G/O)&=0,\qquad
 \pi_{12}(G/O)\cong\mathbb Z,\qquad
 \pi_{13}(G/O)\cong\mathbb Z/3,\notag
\end{align}
and
\begin{equation}\label{eq:eta-action}
 \eta_*a_9=0,\qquad
 \eta_*b_9=j_*(\eta\mu_9)\ne0.
\end{equation}
Here $\eta\in\pi_1^S$ denotes the stable Hopf element,
and $\eta_*$ denotes its action on the homotopy groups
of $G/O$ by stable composition.
These follow from the low stable stems and their compositions,
together with the homotopy exact sequence of $O\to G\to G/O$;
see \cite{Toda1962,MadsenMilgram1979}.  The particular kernel in
\eqref{eq:eta-action} is also recorded in the proof of
\cite[Lemma~3.31]{BasuKasilingamSarkarCP2025}.

We now prove Lemma \ref{Normal_invariant_forh}:
\begin{proof}[Proof of Lemma \ref{Normal_invariant_forh}]
	Let $h:M'\to M$ be an orientation-preserving $\PL$ homeomorphism.
First we have
\[
\eta^{\mathrm{Diff}}(h)
\in
\mathrm{im}\{[M,\PL/O]\to [M,G/O] \}.
\]
Since \(\PL/O\) is \(6\)-connected,
\[
\eta^{\mathrm{Diff}}(h)
\in
\mathcal F^7[M,G/O].
\]
Moreover, since \(\pi_7(G/O)=0\), and
 \(H^8(M;\pi_8(G/O))=0\),
\[
\mathcal F^7[M,G/O]
=
\mathcal F^8[M,G/O]
=
\mathcal F^9[M,G/O],
\]
and therefore
$
\eta^{\mathrm{Diff}}(h)
\in
\mathcal F^9[M,G/O].$
\end{proof}

Before constructing the self-homotopy equivalence as in Lemma \ref{New_self_homotopy}, we compute:  
\begin{lemma}\label{lem:F9}
$
 \mathcal F^9[M,G/O]_{(2)}\cong\mathbb Z/2,$
$ \mathcal F^{10}[M,G/O]_{(2)}=0.$
The nonzero element $\theta_M$ is detected by the permanent cycle
\[
 \bar x\bar y\otimes a_9
 \in E_2^{9,-9}
 =H^9\bigl(M;\pi_9(G/O)_{(2)}\bigr).
\]
\end{lemma}

\begin{proof}
Write $Sq^2\bar y=c\bar x\bar y$, with $c\in\mathbb Z/2$.
Since $M$ is spin, its second Wu class is zero.  The Wu formula and
Cartan formula give
\[
 c=\bigl\langle Sq^2(\bar x^2\bar y),[M]_2\bigr\rangle
   =\bigl\langle v_2(M)\bar x^2\bar y,[M]_2\bigr\rangle=0.
\]

Consequently
\begin{equation}\label{eq:sq2xy}
 Sq^2\bar y=0,\qquad
 Sq^2(\bar x\bar y)=\bar x^2\bar y.
\end{equation}

In the cellular description of $d_2$, the stable attaching component
between cells two dimensions apart is detected by $Sq^2$; its action
on the coefficient group is composition with the stable Hopf map
$\eta$.  This is the usual two-cell description of the cellular
AHSS; compare \cite[Chapter~16]{MosherTangora1968}.
Thus the differential
\[
 d_2:E_2^{9,-9}\longrightarrow E_2^{11,-10}
\]
is, under the generators in \eqref{eq:sq2xy}, the homomorphism
$\eta_*$ in \eqref{eq:eta-action}.  Its kernel is
\[
 \mathbb Z/2\{\bar x\bar y\otimes a_9\}.
\]
The differential entering from $E_2^{7,-8}$ is zero because
$Sq^2\bar y=0$.

For completeness, the possible higher incoming differentials have
sources
\[
 E_r^{9-r,r-10},\qquad 3\le r\le9.
\]
For $r=3,5,7,9$ their coefficient groups are, respectively,
$\pi_7(G/O)$, $\pi_5(G/O)$, $\pi_3(G/O)$, and $\pi_1(G/O)$,
which vanish.  For $r=4,6,8$ their cohomological degrees are
$5,3,1$, where the cohomology of $M$ vanishes.

The outgoing $d_3$ has zero target.  The outgoing $d_4$ has target
\[
 E_4^{13,-12}\cong\mathbb Z_{(2)}:
\]
neither $d_2$ nor $d_3$ has a nonzero source entering this term.
Hence $d_4$ vanishes on the remaining element of order two.
All later outgoing differentials exceed the dimension of $M$.
It follows that
\[
 E_\infty^{9,-9}
 =\mathbb Z/2\{\bar x\bar y\otimes a_9\}.
\]

Finally, it follows that 
\[
\mathcal F^{10}[M,G/O]_{(2)}=0,
\qquad
\mathcal F^{9}[M,G/O]_{(2)}
=E_\infty^{9,-9}\cong\Z/2.
\]
by
$H^{10}(M;-)=H^{12}(M;-)=0$ and 
$\pi_{11}(G/O)_{(2)}=\pi_{13}(G/O)_{(2)}=0.$
\end{proof}

\subsubsection{The pinch map and its Pontryagin--Thom class}
$\quad$

In this subsubsection, we prove Lemma \ref{New_self_homotopy}.
Let $i:\mathbb{CP}^3\to M$ be the projective-stage map constructed
earlier, and put $x_0=i^*x$.  Since
$
2\dim\mathbb{CP}^3<\dim M,$
general position allows us to take \(i\) as a smooth embedding.
Let
\[
 H:S^7\longrightarrow\mathbb{CP}^3
\]
be the complex Hopf projection, which is also the attaching map of
the top cell of $\mathbb{CP}^4$.  Choose an unstable representative
\[
 \nu_7^2:S^{13}\longrightarrow S^7
\]
of $\nu^2\in\pi_6^S$.

Choose a disk in $M$ disjoint from $i(\mathbb{CP}^3)$ and let
$c:M\to M\vee S^{13}$ be the corresponding pinch map.  Define
\begin{equation}\label{eq:pinch}
 \sigma=i\circ H\circ\nu_7^2,
 \qquad
 f_{\nu^2}^M=(\operatorname{id}_M\vee\sigma)\circ c.
\end{equation}

We keep the rank normalization of Thom spectra explicit.  Set
\[
 \nu_M^0=[\varepsilon^{13}]-[TM]
 \in\widetilde{KO}(M).
\]
This is the virtual rank-zero stable normal class.  Atiyah duality
\cite[Theorem~3.3]{Atiyah1961} gives
\begin{equation}\label{eq:rank-zero-duality}
 M^{\nu_M^0}\simeq\Sigma^{13}D(M_+)
\end{equation}
where \(M^{\nu_M^0}\) is the Thom spectrum of \(\nu_M^0\) and
\(D(M_+)\) is the Spanier--Whitehead dual of \(M_+=M\sqcup\{*\}\).

All homotopy groups of Thom spectra and their cell quotients below
are stable; in such expressions $S^d$ denotes the sphere spectrum.
Thus $\pi_{13}(S^4)=\pi_9^S$ is not an assertion about the
unstable homotopy group of the four-sphere.

Let $L$ be the complex line bundle over $\CP^3$ with $c_1(L)=x_0$, and put
\[
 u=[L_{\mathbb R}]-2.
\]
Fujii's calculation gives
\[
\widetilde{KO}(\mathbb{CP}^3)=\mathbb Z\{u\},
\qquad
p_1(u)=x_0^2.
\]
Recall that $\lambda$ denotes the first Pontryagin coefficient of $M$. Therefore
\[
i^*\nu_M^0=-\lambda u.
\]
Extend this class to $\mathbb{CP}^4$ by the same expression
$\xi=-\lambda u$, and let
\[
 T_j=(\mathbb{CP}^j)^{\xi|_{\mathbb{CP}^j}},
 \qquad 0\le j\le4
\]
denote the Thom spectrum of $\xi|_{\mathbb{CP}^j}$.
Then we have a Thom cofibration
\begin{equation}\label{eq:Thom-attachment}
 S^7\xrightarrow{A_\lambda}T_3\longrightarrow T_4.
\end{equation}
Fix a stable bundle identification
\[
 \xi|_{\mathbb{CP}^3}\cong i^*\nu_M^0,
\]
and denote the induced Thom map by
\[
 i^\nu:T_3\longrightarrow M^{\nu_M^0}.
\]

\begin{lemma}\label{lem:pinch-PT}
The map $f_{\nu^2}^M$ is an orientation-preserving tangential
self-homotopy equivalence.  For a compatible choice of tangential
normal data, its reduced tangential normal invariant is the
$S$-dual of
\begin{equation}\label{eq:Anu2}
 \rho_\sigma=i^\nu_*(A_\lambda\circ\nu^2)
 \in\pi_{13}(M^{\nu_M^0}).
\end{equation}
Here $\nu^2$ in a composite of spectra denotes the stabilization
of $\nu_7^2$.
\end{lemma}

\begin{proof}
The map $\sigma$ factors through a six-dimensional complex, and
therefore induces zero on $H_{13}$.  The pinch map consequently
induces the identity on integral homology.  Since $M$ is simply
connected, it is a degree-one homotopy equivalence.

For $v\in\widetilde{KO}(M)$ the pinch coaction gives
\[
 (f_{\nu^2}^M)^*v=v+q^*(\sigma^*v),
\]
where $q:M\to S^{13}$ collapses the complement of the chosen disk.
Bott periodicity gives
\[
 \widetilde{KO}(S^{13})=0.
\]
Hence the second term vanishes, and in particular
\[
 (f_{\nu^2}^M)^*TM\cong_{\mathrm{st}}TM.
\]
Let \(b\) denote a stable bundle isomorphism covering
\(f_{\nu^2}^M\).  By \cite[Lemma~7.4]{CrowleyHambleton2015},
\[
 \eta^t(f_{\nu^2}^M,b)=1\ast  \widehat D(\rho_\sigma)
 \in [M,SG].
\]
Equation~(7.2) of \cite{CrowleyHambleton2015} gives
$$
 \rho_\sigma=i^\nu_*(\rho_y),
 \qquad
 y=H\circ\nu_7^2,
$$
where \(\rho_y\) is the Pontryagin--Thom class determined by \(y\)
together with the chosen trivialization of
\(y^*(\xi|_{\CP^3})\).

We now identify this class explicitly.  Recall that
\(H:S^7\to\CP^3\) is the attaching map of the top cell of
\(\CP^4\), and that \(\xi=-\lambda u\) has been extended from
\(\CP^3\) to \(\CP^4\).  Let
$$
 \chi:(D^8,S^7)\longrightarrow(\CP^4,\CP^3)
$$
be a characteristic map for the top cell, so that
\(\chi|_{S^7}=H\).  Since \(D^8\) is contractible,
\(\chi^*\xi\) is stably trivial.  Restricting such a trivialization
to \(S^7\) gives a preferred stable trivialization of
\(H^*(\xi|_{\CP^3})\).  With this choice, the induced Thom map
$$
 (S^7)^{H^*(\xi|_{\CP^3})}\longrightarrow
 {(\CP^3)}^{\xi|_{\CP^3}}=T_3
$$
is precisely the Thom attaching map
$$
 A_\lambda:S^7\longrightarrow T_3
$$
appearing in \eqref{eq:Thom-attachment}.
Pulling this trivialization back along
\(\nu_7^2:S^{13}\to S^7\), the Pontryagin--Thom construction is
natural under composition.  Hence
$
 \rho_y=A_\lambda\circ\nu^2
 \in\pi_{13}(T_3).
$
\end{proof}

\begin{proposition}\label{prop:pinch}
$
 \eta^{\mathrm{Diff}}(f_{\nu^2}^M)_{(2)}
 =\theta_M\ne0.$
More precisely,
\[
 \operatorname{gr}_9
 \eta^{\mathrm{Diff}}(f_{\nu^2}^M)_{(2)}
 =
 \bar x\bar y\otimes a_9\in E_\infty^{9,-9}.
\]

\end{proposition}

\begin{proof}
We first compute the relevant component of the Thom attaching map.
Let $U$ be the mod-$2$ Thom class of the virtual rank-zero bundle
$\xi$ over $\mathbb{CP}^4$.  Since
\[
 w(\xi)=(1+\bar x_0)^{-\lambda}
\]
and $4\mid\lambda$, we have
\begin{equation}\label{eq:w2w4}
 w_2(\xi)=w_4(\xi)=0.
\end{equation}
The Thom identity and the Cartan formula therefore give
\begin{align}
 Sq^2(\bar x_0^2U)&=0,
 \label{eq:Thom-Sq2-bottom}\\
 Sq^2(\bar x_0^3U)&=\bar x_0^4U,\notag\\
 Sq^4(\bar x_0^2U)&=\bar x_0^4U.\notag
\end{align}

The quotient $T_3/T_1$ has cells in dimensions $4$ and $6$.
The first identity in \eqref{eq:Thom-Sq2-bottom} shows that the
stable attaching map of the $6$-cell to the $4$-cell is trivial,
since
\[
 \pi_1^S\cong\mathbb Z/2\{\eta\}
\]
and $Sq^2$ detects this attaching map.  Hence, stably and
$2$-locally,
\begin{equation}\label{eq:wedge}
 T_3/T_1\simeq S^4\vee S^6.
\end{equation}
With respect to this splitting, the image of the Thom attaching map
$A_\lambda:S^7\to T_3$ has the form
\begin{equation}\label{eq:Acomponents}
 S^7\xrightarrow{A_\lambda}T_3
 \longrightarrow T_3/T_1\simeq S^4\vee S^6,
 \qquad
 A_\lambda=(b_\lambda\nu,\eta),
\end{equation}
where \(b_\lambda\in\mathbb Z/8\) after \(2\)-localization.  The
second component is nonzero because
\[
Sq^2(\bar x_0^3U)=\bar x_0^4U.
\]
Likewise,
\[
Sq^4(\bar x_0^2U)=\bar x_0^4U
\]
shows that
\begin{equation}\label{eq:bodd}
 b_\lambda\equiv1\pmod2.
\end{equation}

Put
\[
\beta=A_\lambda\circ\nu^2
\in\pi_{13}(T_3)_{(2)}.
\]
Using
\[
 \eta\nu^2=0,\qquad
 \nu\nu^2=\nu^3,\qquad
 2\nu^3=0,\qquad
 \nu^3\ne0
\]
from \cite{Toda1962}, equations
\eqref{eq:Acomponents} and \eqref{eq:bodd} give
\[
 \beta\longmapsto(\nu^3,0)
 \in\pi_{13}(S^4\vee S^6)_{(2)}.
\]
In particular, the image of $\beta$ in
\[
 \pi_{13}(T_3/T_2)_{(2)}
 \cong\pi_{13}(S^6)_{(2)}
\]
vanishes, so \(\beta\) lifts to a class
\[
\widetilde\beta\in\pi_{13}(T_2)_{(2)}.
\]
Its image under
\[
T_2\longrightarrow T_2/T_1\simeq S^4
\]
is
\begin{equation}\label{eq:nu3-fourcell}
 \nu^3\in\pi_{13}(S^4)_{(2)}
       =\pi_9^S{}_{(2)}.
\end{equation}

The choice of the lift \(\widetilde\beta\) is immaterial for this
projection.  Indeed, the ambiguity comes from
\[
\pi_{14}(T_3/T_2)_{(2)}
\cong
\pi_{14}(S^6)_{(2)},
\]
and its image in \(\pi_{13}(T_2)\), after projection to
\(T_2/T_1\simeq S^4\), factors through the attaching map of the
\(6\)-cell to the \(4\)-cell.  This attaching map is stably trivial by
\eqref{eq:wedge}.  Hence every lift has the same nonzero image
\(\nu^3\) in \(\pi_{13}(S^4)_{(2)}\).
Consider the homological AHSS
$$E^2_{p,q}=H_p(T_3;\pi_q^s)\Longrightarrow \pi_{p+q}(T_3)$$
associated with the filtrations
$$\mathcal F_{p}\pi_{13}(T_3)=\im(\pi_{13}(T^{(p)}_3)\to \pi_{13}(T_3)).$$
It follows that \(\beta\) does not lie in the filtration coming from
\(T_1=T^{(2)}_3\).  Thus \(\beta\) has filtration exactly \(4\).  If $e_4$ denotes the cellular generator of
$H_4(\mathbb{CP}^3;\mathbb Z/2)$ corresponding to the $4$-cell, then
\begin{equation}\label{eq:beta-gr4}
 \operatorname{gr}_4(\beta)
 =
 e_4\otimes\nu^3
 \ne0
 \in E^\infty_{4,9}(T_3)_{(2)}.
\end{equation}

Recall the projective-stage embedding
 $i:\mathbb{CP}^3\hookrightarrow M$
  and let $\nu_i$ be its rank-$7$
normal bundle.  Atiyah duality gives
\[
 \Sigma^{13}D(T_3)\simeq\operatorname{Th}(\nu_i).
\]
where $D(T_3)$ is the Spanier--Whitehead dual of $T_3$, $\operatorname{Th}(\nu_i)$
denotes the Thom space of $\nu_i$.

Let
$\widetilde\pi_S^0(X)
:=
[\Sigma^\infty X,\mathbb S]$
denote reduced stable cohomotopy, where \(\mathbb S\) is the sphere
spectrum.  By the spectrum--infinite-loop-space adjunction, for a
connected based finite CW complex \(X\),
\[
\widetilde\pi_S^0(X)
\cong
[X,(QS^0)_0].
\]
Thus
\[
\alpha_\beta :=\widehat D(\beta)
\in
[\operatorname{Th}(\nu_i),(QS^0)_0]_*
\cong
\widetilde\pi_S^0(\operatorname{Th}(\nu_i)).
\]

Spanier--Whitehead duality identifies the homological and
cohomological Atiyah--Hirzebruch spectral sequences associated to
dual cellular filtrations; see
\cite[Theorem~A and Definition~3.2]{Hans2026}.  After the
$13$-dimensional shift above, the bidegrees correspond by
\[
 (p,q)\longleftrightarrow(13-p,-q).
\]

In the cohomological AHSS for $\widetilde\pi_S^0(\operatorname{Th}(\nu_i))$,
 \eqref{eq:beta-gr4} gives
\begin{equation}\label{eq:alpha-beta-gr9}
 \operatorname{gr}_9(\alpha_\beta)
 =
 \bar x_0\,\bar U_i\otimes\nu^3
 \ne0
 \qquad\text{in }
 E_\infty^{9,-9}
 \bigl(\operatorname{Th}(\nu_i);\mathbb S\bigr)_{(2)},
\end{equation}
where
$
 U_i\in
 H^7(\operatorname{Th}(\nu_i);\mathbb Z)$
is the Thom class and bars denote mod-$2$ reduction.

Recall the Thom map
\[
i^\nu:T_3\longrightarrow M^{\nu_M^0}
\]
induced by the inclusion \(i\).  Its Spanier--Whitehead dual is the
Pontryagin--Thom collapse map
\[
 i^!:M_+\longrightarrow\operatorname{Th}(\nu_i);
\]
see \cite[Lemma~7.3(i)]{CrowleyHambleton2015}.  By the
Pontryagin--Thom identification discussed above,
\[
 \rho_\sigma=i^\nu_*(\beta),
 \qquad
 \sigma=i\circ H\circ\nu_7^2.
\]
Naturality of Spanier--Whitehead duality therefore gives
\begin{equation}\label{eq:dual-rho-sigma}
 \widehat D(\rho_\sigma)
 =
 (i^!)^*\alpha_\beta.
\end{equation}

Recall the Gysin homomorphism
\begin{equation}\label{eq:collapse-gysin}
i_!(b):=(i^!)^*(bU_i):
H^r(\CP^3;\Z)\longrightarrow H^{r+7}(M;\Z).
\end{equation}
It satisfies
\[
\bigl\langle a\cup i_!(b),[M]\bigr\rangle
=
\bigl\langle i^*a\cup b,[\CP^3]\bigr\rangle.
\]
Put \(x_0=i^*x\).  Choose the sign of the Thom class \(U_i\) so that
\[
i_!(1)=y.
\]
The projection formula then gives
\[
i_!(x_0)
=
i_!(i^*x)
=
x\,i_!(1)
=
xy.
\]
Consequently,
\begin{equation}\label{eq:gysin-x0}
 (i^!)^*(\bar x_0\,\bar U_i)
 =
 \bar x\bar y.
\end{equation}

We now return to the normal invariant.  By
\cite[Lemmas~7.3--7.4]{CrowleyHambleton2015} and Lemma \ref{lem:pinch-PT}, the tangential normal
invariant of the pinch map is
\[
 \eta^t(f_{\nu^2}^M,b)=[1]*\widehat D(\rho_\sigma)
 =
 [1]*(i^!)^*\alpha_\beta
 \in[M,SG].
\]

Consider the following commutative diagram
\begin{equation}
\begin{tikzcd}[column sep=large,row sep=large]
 {[\mathrm{Th}(\nu_i), SG]}
\arrow[r,"(i^!)^\ast"]
\arrow[d,"j_\ast"]
&
{[M, SG]}\arrow[d,"j_\ast"]
\\
{[\mathrm{Th}(\nu_i),G/O]}
\arrow[r,"(i^!)^\ast"]
&
{[M, G/O]}. 
\end{tikzcd}
\end{equation}
By the naturality of AHSS, \eqref{eq:alpha-beta-gr9} gives 
$$\eta^{\Diff}(f_{\nu^2}^M)=j_\ast([1]*\widehat D(\rho_\sigma))\in \mathcal F^9[M,G/O].$$
  By \eqref{eq:alpha-beta-gr9}, \eqref{eq:gysin-x0} and
$j_*(\nu^3)=a_9$, it follows that $\eta^{\Diff}(f_{\nu^2}^M)$ is detected by
\[
 \bar x\bar y\otimes a_9\in E_2^{9,-9}.
\]
By \cref{lem:F9}, 
$\eta^{\mathrm{Diff}}(f_{\nu^2}^M)_{(2)}
=
\theta_M\ne0.$
\end{proof}

This completes the proof for Lemma \ref{New_self_homotopy}.

\subsubsection{Proof for Lemma \ref{3-primary_normal_in}}

Assume that there exists a homotopy equivalence
$g:M'\to M$ such that $\eta^{\Diff}(g)\in \mathcal F^9[M,G/O]$ and $\eta^{\Diff}(g)_{(2)}=0$.
Above the eight-skeleton the only nonzero cohomological degrees of
$M$ are $9,11,13$.
The coefficient groups in degree eleven vanish,
the degree-nine group is $\Z/2\oplus \Z/2$, and
\[
 \pi_{13}(G/O)\cong\mathbb Z/3.
\]
The AHSS therefore implies
\[
 \zeta:=\eta^{\Diff}(g)\in\mathcal F^{13}[M,G/O].
\]
For the degree-one top-cell collapse $q:M\to S^{13}$, exactness
of the skeletal cofiber sequence gives
\begin{equation}\label{eq:top-filtration-image}
 \mathcal F^{13}[M,G/O]
 =
 \operatorname{im}
 \bigl(q^*:\pi_{13}(G/O)\to[M,G/O]\bigr).
\end{equation}

The smooth surgery sequence for the sphere contains
\[
L_{14}(e)\longrightarrow\Theta_{13}
\xrightarrow{\eta_S}\pi_{13}(G/O)
\longrightarrow L_{13}(e),
\]
and \(\eta_S\) is an isomorphism; see
\cite{KervaireMilnor1963}.

Choose $\Sigma\in\Theta_{13}$ so that
$ q^*\eta_S(\Sigma)=-\zeta.$
Connected sum produces a homotopy equivalence
$
 g_2:M'\#\Sigma\longrightarrow M
$
with
\[
 \eta^{\mathrm{Diff}}(g_2)
 =
 \eta^{\mathrm{Diff}}(g)
 +q^*\eta_S(\Sigma)
 =0.
\]
Exactness puts $(M'\#\Sigma,g_2)$ in the orbit of the base point under the
$L_{14}(e)$-action. This action is connected sum with an element of $bP_{14}$;
see \cite{Wall1999,KervaireMilnor1963}. Here $bP_{14}=0$.
Thus $M'\#\Sigma\cong_{\Diff}M$.

\begin{remark}
For a general even $\lambda$ the same Thom calculation gives
\[
 w_4(-\lambda u)
 =
 \left(\frac{\lambda}{2}\bmod2\right)\bar x_0^2,
\]
and
\[
 \operatorname{gr}_9
 \eta^{\mathrm{Diff}}(f_{\nu^2}^M)_{(2)}
 =
 \left(1+\frac{\lambda}{2}\bmod2\right)
 \bar x\bar y\otimes a_9.
\]
Thus this particular pinch removes the filtration-nine obstruction
when $\lambda\equiv0\pmod4$.
Its failure to detect that obstruction when
$\lambda\equiv2\pmod4$.
\end{remark}

\section{The smooth classification in dimension nine}\label{sec:9smooth}

Throughout this section 
$
 M=M_{a,b},$
 $\lambda=\sum_{i=1}^{3}(a_i^2+b_i^2).$
By \cref{prop:tangent},
\begin{equation}\label{eq:p1lambda9}
 p_1(M)=\lambda x^2,
 \qquad
 w_2(M)\equiv\lambda x\pmod2.
\end{equation}
Thus $M$ is non-spin exactly when $\lambda$ is odd.

\subsection{\texorpdfstring{The top cell and stable real $K$-theory}{The top cell and stable real K-theory}}\label{subsec:9-cell}

We first record a convenient cellular model.

\begin{lemma}\label{lem:9d-cell}
The manifold \(M\) has a CW model
\[
M\simeq(\CP^2\vee S^5)\cup_{\beta}e^7\cup_{\alpha}e^9.
\]
If \(K=M^{(7)}\), then
\[
K\simeq(\CP^2\vee S^5)\cup_{\beta}e^7.
\]
\end{lemma}

\begin{proof}
This is a standard minimal decomposition for a simply connected complex.
\end{proof}

Let \(\bar x,\bar y\) denote the mod-\(2\) reductions of \(x,y\).

\begin{lemma}\label{lem:9d-eta-topcell}
Assume that $M$ is non-spin.  Let
$c\colon K\to K/K^{(6)}\simeq S^7$ collapse the $6$-skeleton.  Then the
homotopy class of
\[
 S^8\xrightarrow{\alpha}K\xrightarrow{c}S^7
\]
is the nonzero element
$\eta_7\in\pi_8(S^7)\cong\Z/2$.

\end{lemma}

\begin{proof}
Since $M$ is non-spin, \eqref{eq:p1lambda9} gives
$v_2(M)=w_2(M)=\bar x$.  The Wu formula gives
\[
 \langle Sq^2(\bar x\bar y),[M]_2\rangle
 =\langle v_2(M)\bar x\bar y,[M]_2\rangle
 =\langle \bar x^2\bar y,[M]_2\rangle=1.
\]
Thus
\begin{equation}\label{eq:Sq2xy}
Sq^2(\bar x\bar y)=\bar x^2\bar y\ne0.
\end{equation}
The corresponding functional Steenrod square detects the projection
of the top-cell attaching map to the \(7\)-cell.  Hence
\(c\circ\alpha\) is nonzero in
\(\pi_8(S^7)\cong\mathbb Z/2\); therefore
\[
c\circ\alpha=\eta_7.
\]
See, for example, \cite[Chapter~16]{MosherTangora1968} for this
two-cell detection of the attaching map.
\end{proof}

\subsection{Vanishing of the smooth normal invariant}\label{subsec:9-normal}

Suppose that \(\lambda=\lambda'\).  By \cref{thm:main}, there
is a PL homeomorphism
\[
h:M'\longrightarrow M.
\]
Choose an orientation on \(M\), and give \(M'\) the pulled-back
orientation, so that \(h\) is orientation-preserving.  

\begin{proposition}\label{prop:9d-normal-invariant-zero}
Let $M$ and $M'$ be non-spin nine-dimensional orbit manifolds with
$\lambda=\lambda'$.  For a $\PL$ homeomorphism
$h\colon M'\to M$,
$
 \eta^{\mathrm{Diff}}(h)=0\in[M,G/O].$
\end{proposition}

\begin{proof}
The normal invariant of a PL homeomorphism lies in the image of the canonical
map
\[
 [M,\PL/O]\longrightarrow[M,G/O].
\]
Let $i\colon K=M^{(7)}\to M$ denote the inclusion.  The space $\PL/O$ is
$6$-connected, so a map $K\to \PL/O$ is detected in dimension $7$.  Since
$\pi_7(G/O)=0$, the induced homomorphism
\[
 [K,\PL/O]\longrightarrow[K,G/O]
\]
is zero.  Therefore
\begin{equation}\label{eq:eta-restrict-zero}
 i^*\eta^{\mathrm{Diff}}(h)=0.
\end{equation}

Apply $[-,G/O]$ to the cofibration
\[
 K\xrightarrow{i} M\xrightarrow{q}S^9\xrightarrow{\Sigma \alpha}\Sigma K.
\]
We obtain the exact segment
\begin{equation}\label{eq:GO-puppe}
 [\Sigma K,G/O]\xrightarrow{(\Sigma\alpha)^*}
 \pi_9(G/O)\xrightarrow{q^*}[M,G/O]\xrightarrow{i^*}[K,G/O].
\end{equation}
By \cref{lem:9d-eta-topcell}, the composite
$\Sigma c\circ\Sigma\alpha$ is $\eta_8$.  Basu--Kasilingam--Sarkar prove that
\begin{equation}\label{eq:eta-GO-surj}
 \eta^*\colon\pi_8(G/O)\longrightarrow\pi_9(G/O)
\end{equation}
is onto \cite[Lemma~3.7]{BasuKasilingamSarkar2026}.  Hence \((\Sigma\alpha)^*\) in \eqref{eq:GO-puppe} is surjective.
Indeed, for \(z\in\pi_9(G/O)\), choose
\(\gamma\in\pi_8(G/O)\) with \(\gamma\circ\eta_8=z\).  The class
represented by
\[
\gamma\circ\Sigma c:\Sigma K\longrightarrow G/O
\]
then maps to \(z\) under \((\Sigma\alpha)^*\).  Thus \(q^*=0\).
Together with \eqref{eq:eta-restrict-zero} and exactness of
\eqref{eq:GO-puppe}, this gives
\(\eta^{\mathrm{Diff}}(h)=0\).
\end{proof}

\subsection{\texorpdfstring{The final $L_{10}$-obstruction}{The final L10-obstruction}}\label{subsec:9-L10}

The relevant part of the smooth surgery exact sequence is
\begin{equation}\label{eq:smooth-surgery-9d}
 [\Sigma M,G/O]\xrightarrow{\Theta_M^{\mathrm{Diff}}}
 L_{10}(e)\longrightarrow\mathcal S^{\mathrm{Diff}}(M)
 \xrightarrow{\eta^{\mathrm{Diff}}}[M,G/O]\longrightarrow L_9(e),
\end{equation}
where $L_{10}(e)\cong\Z/2$ and $L_9(e)=0$.  We now show that the first map in
\eqref{eq:smooth-surgery-9d} is onto.  The construction follows the proof of
\cite[Theorem~5.7]{BasuKasilingamSarkar2026}, but the required class in
$[\Sigma M,G/O]$ has to be built from the cell structure of $M$ rather than
from a product decomposition.

\begin{lemma}\label{lem:9d-GO-class}
Assume that $M$ is non-spin, and let
$k_6\in H^6(G/O;\Z/2)$ denote the smooth Kervaire class.
There exists $u\in[\Sigma M,G/O]$ such that
\[
 u^*(k_6)=\sigma\bar y\ne0
 \quad\text{in }H^6(\Sigma M;\Z/2),
\]
where $\sigma\bar y$ is the suspension of the degree-$5$ generator.
\end{lemma}

\begin{proof}
On
$
\Sigma(\CP^2\vee S^5)=\Sigma\CP^2\vee S^6$,
define \(u_0\) to be null on \(\Sigma\CP^2\) and to represent
\(\nu^2\in\pi_6(G/O)\) on \(S^6\).  Since
\(\pi_7(G/O)=0\), the map \(u_0\) extends across the \(8\)-cell to
\[
u_1:\Sigma K\longrightarrow G/O.
\]

Now
$\Sigma M=\Sigma K\cup_{\Sigma\alpha}e^{10}.$
If \(u_1\circ\Sigma\alpha\) is null-homotopic, \(u_1\) already extends
to the required map \(u:\Sigma M\to G/O\).  Otherwise, write
\[
z=[u_1\circ\Sigma\alpha]\in\pi_9(G/O).
\]
By the surjectivity of \eqref{eq:eta-GO-surj}, choose
\(\gamma\in\pi_8(G/O)\) such that
\[
\gamma\circ\eta_8=-z.
\]
Using the \(H\)-space structure on \(G/O\), define
\[
u'_1
=
u_1+(\gamma\circ\Sigma c):
\Sigma K\longrightarrow G/O.
\]
Since \(\Sigma c\circ\Sigma\alpha=\eta_8\), we have
\[
u'_1\circ\Sigma\alpha\simeq *,
\]
and hence \(u'_1\) extends to a map
\(u:\Sigma M\to G/O\).

The correction term \(\gamma\circ\Sigma c\) factors through \(S^8\)
and therefore has trivial pullback in degree \(6\).  It follows that
the degree-\(6\) Kervaire class is unchanged, i.e.,
$u^*(k_6)=u_0^*(k_6)=\sigma\bar y\ne0.$
\end{proof}

\begin{proposition}\label{prop:9d-theta-surjective}
For every non-spin nine-dimensional orbit manifold $M$, the surgery
obstruction homomorphism
\[
\Theta_M^{\mathrm{Diff}}:
[\Sigma M,G/O]\longrightarrow L_{10}(e)\cong\Z/2
\]
is surjective.
\end{proposition}

\begin{proof}
Let \(u\in[\Sigma M,G/O]\) be the class from
\cref{lem:9d-GO-class}, so that
$
u^*(k_6)=\sigma\bar y\ne0.$
Let
\[
p:M\times S^1\longrightarrow M\wedge S^1=\Sigma M
\]
be the quotient map and put \(F=u\circ p\).  Then
\[
F^*(k_6)=\bar y\,t,
\]
where \(t\in H^1(S^1;\Z/2)\) is the generator.
By the commutative diagram
\[
\begin{tikzcd}[column sep=large,row sep=large]
 {[\Sigma M, G/O]}
\arrow[r,"p^\ast"]
\arrow[d,"\Theta_M^{\mathrm{Diff}}"]
&
{[M\times S^1, G/O]}\arrow[d,"\Theta^{\mathrm{Diff}}"]
\\
L_{10}(e)
\arrow[r,"\cong"]
&
L_{10}(\Z), 
\end{tikzcd}
\]
 $0\ne \Theta_M^{\mathrm{Diff}}(u)\in L_{10}(e)$ is equivalent to $0\ne \Theta^{\mathrm{Diff}}(F)\in L_{10}(\Z)$.

The smooth surgery obstruction may be evaluated by the Kervaire-class
formula \cite{RourSull}:
\begin{equation}\label{eq:Kervaire-obstruction-9d}
 \Theta^{\mathrm{Diff}}(F)=
 \left\langle V(M\times S^1)^2\,
 F^*\!\left(\sum_{j\ge0}k_{2^{j+1}-2}\right),
 [M\times S^1]_2\right\rangle
\end{equation}
where $V$ is the total Wu class, $k_{2^{j+1}-2}\in H^{2^{j+1}-2}(G/O;\Z/2)$.

Since $H^2(\Sigma M;\Z/2)=H^1(M;\Z/2)=0$, the degree-two
Kervaire class has zero pullback.  Classes of degree at least $14$ cannot
contribute in dimension $10$.  Thus the only contributing term is
\[
\Theta^{\mathrm{Diff}}(F)
=
\left\langle
(V(M\times S^1)^2)_4\,F^*(k_6),
[M\times S^1]_2
\right\rangle.
\]
Since \(M\) is non-spin,
\[
v_2(M)=w_2(M)=\bar x,
\]
and hence
\[
(V(M\times S^1)^2)_4=v_2(M)^2=\bar x^2.
\]
Therefore
\[
\left\langle
\bar x^2\bar y\,t,[M\times S^1]_2
\right\rangle
=
\left\langle\bar x^2\bar y,[M]_2\right\rangle
\left\langle t,[S^1]_2\right\rangle
=1.
\]
This gives that $\Theta_M^{\mathrm{Diff}}$ is surjective.
\end{proof}

\begin{theorem}\label{thm:9d-smooth}
Let \(M_{a,b}\) and \(M_{a',b'}\) be orbit manifolds of free diagonal
linear \(S^1\)-actions on \(S^5\times S^5\).  Suppose that both
\(\lambda(a,b)\) and \(\lambda(a',b')\) are odd.  Then
\[
 M_{a,b}\cong_{\mathrm{Diff}}M_{a',b'}
 \quad\Longleftrightarrow\quad
 \lambda(a,b)=\lambda(a',b').
\]
\end{theorem}

\begin{proof}
Necessity follows from \eqref{eq:p1lambda9}, because $p_1$ is preserved by
diffeomorphisms.  Conversely, suppose
$\lambda(a,b)=\lambda(a',b')$ is odd.  By \cref{thm:main}, there is a PL
homeomorphism $h\colon M_{a',b'}\to M_{a,b}$.  By
\cref{prop:9d-normal-invariant-zero},
$\eta^{\mathrm{Diff}}(h)=0$.
By \cref{prop:9d-theta-surjective}, the preceding homomorphism
\[
[\Sigma M,G/O]\longrightarrow L_{10}(e)
\]
is surjective.  Hence the map
\(L_{10}(e)\to\mathcal S^{\mathrm{Diff}}(M)\) is trivial.
Exactness of \eqref{eq:smooth-surgery-9d} now shows that
\((M_{a',b'},h)\) represents the base point of the smooth structure
set of \(M_{a,b}\).  Thus the two manifolds are diffeomorphic.
\end{proof}

\section{Sphere-bundle realization and homotopy classification}\label{sec:sphere-realization}

The projective stage gives a geometric model for every orbit manifold.
It also makes the homotopy classification transparent: after passing to
the associated sphere bundle, only the stable \(J\)-class of the bundle
over \(\CP^n\) remains.

\begin{theorem}\label{thm:sphere-realization}
Let \(n\ge2\), and let \(M=M_{a,b}\) be an orbit manifold of a free
diagonal linear \(S^1\)-action on
\(S^{2n+1}\times S^{2n+1}\).  There is an oriented real vector bundle
\[
 \eta_{a,b}\longrightarrow\CP^n
\]
of rank \(2n+2\) such that
\[
 M\cong_{\PL}S(\eta_{a,b}).
\]
If \(i:\CP^n\to M\) is the projective-stage map, then
\begin{equation}\label{eq:eta-stable-class}
 [T\CP^n]+[\eta_{a,b}]=[i^*TM]+[\eps^1]
 \qquad\text{in }KO(\CP^n),
\end{equation}
and hence
\begin{equation}\label{eq:eta-pontryagin}
 p(\eta_{a,b})
 =\frac{i^*p(TM)}{p(T\CP^n)}
 =\frac{\displaystyle
 \prod_{i=1}^{n+1}(1+a_i^2x^2)
 \prod_{j=1}^{n+1}(1+b_j^2x^2)}
 {(1+x^2)^{n+1}}.
\end{equation}
\end{theorem}

\begin{proof}
The virtual bundle
\[
 [i^*TM]+[\eps^1]-[T\CP^n]
\]
has rank \(2n+2\).  Since \(\dim\CP^n=2n\), the stable range for real
vector bundles represents it by an actual oriented rank-\((2n+2)\)
bundle \(\eta_{a,b}\); see \cite[Chapter~15]{Switzer1975}.  This gives
\eqref{eq:eta-stable-class}, and \eqref{eq:eta-pontryagin} follows
immediately from \cref{prop:tangent} and
\(p(T\CP^n)=(1+x^2)^{n+1}\).

Set \(N=S(\eta_{a,b})\) and let \(\pi:N\to\CP^n\) be the bundle
projection.  Since \(\operatorname{rank}\eta_{a,b}=2n+2>\dim\CP^n\),
the bundle admits a nowhere-zero section \(s:\CP^n\to N\).  Hence
\[
 H^*(N;\Z)\cong
 \Z[x_N,y_N]/(x_N^{n+1},y_N^2),
 \qquad |x_N|=2,\quad |y_N|=2n+1,
\]
where \(x_N=\pi^*x\).  The standard tangent-bundle identity for a
sphere bundle is
\begin{equation}\label{eq:sphere-tangent}
 TN\oplus\eps^1
 \cong_{\mathrm{st}}
 \pi^*(T\CP^n\oplus\eta_{a,b}).
\end{equation}
Using \eqref{eq:eta-stable-class}, we obtain
\[
 TN\oplus\eps^1
 \cong_{\mathrm{st}}
 \pi^*(i^*TM\oplus\eps^1).
\]
Therefore
\begin{equation}\label{eq:sphere-pk}
 p_k(N)=\Pi_k(a,b)\,x_N^{2k}
 \qquad(1\le k\le\lfloor n/2\rfloor).
\end{equation}
Pulling \eqref{eq:sphere-tangent} back along the section also gives
\[
 s^*TN\cong_{\mathrm{st}}i^*TM.
\]
Thus, when \(n=4r+1\), the top \(\Z/2\)-coefficient in
\(KO(\CP^n)\) agrees as well and is equal to \(\Delta_r(a,b)\).

Consequently \(M\) and \(N\) have identical stable normal data on the
projective stage and hence the same normal \(2n\)-type.  The canonical
normal smoothing of \(M\) lies in \(F_{2n,2n+1}\) by
\cref{prop:low-filtration}.  The stable normal bundle of \(N\) is
pulled back from \(\CP^n\), so its normal smoothing factors through the
\(2n\)-skeleton and lies in the same filtration.  By
\cref{prop:PLfiltration-zero}, both classes vanish after passage to PL
bordism.  The classification argument in the proof of
\cref{thm:main} now applies directly: \cref{prop:kreck-surgery}, with
\(CAT=\PL\), changes the resulting normal bordism relative to its boundary
into a PL \(s\)-cobordism.  Hence \(M\cong_{\PL}N\).
\end{proof}

\begin{lemma}[The general homotopy criterion]\label{rem:general-J-homotopy}
Let
\[
 \gamma_{a,b}
 = [j_n^*\Gamma_{a,b,\R}]-[\eps^{4n+4}]
 \in\widetilde{KO}(\CP^n).
\]
Then
\begin{equation}\label{eq:general-J-homotopy}
 M_{a,b}\simeq M_{a',b'}
 \quad\Longleftrightarrow\quad
 J(\gamma_{a,b})=J(\gamma_{a',b'})
 \quad\text{in }\widetilde{J}(\CP^n).
\end{equation}
\end{lemma}
\begin{proof}
If \(\widetilde\eta_{a,b}\) is the reduced stable class of
\(\eta_{a,b}\), then \eqref{eq:eta-stable-class} and
\cref{prop:tangent} give
\[
 \widetilde\eta_{a,b}
 =\gamma_{a,b}-\widetilde{T\CP^n}.
\]
Thus equality in \eqref{eq:general-J-homotopy} is precisely equality of
the stable spherical fibrations associated with the two sphere bundles.
Since the bundles have rank \(2n+2\) over the \(2n\)-complex \(\CP^n\),
Freudenthal stabilization is already in the stable range, so the unit
sphere bundles are fiber-homotopy equivalent.  Conversely, a homotopy
equivalence of the total spaces preserves the Spivak normal spherical
fibration; pulling back to the projective stages, which are unique up
to complex conjugation, gives the reverse implication.  Complex
conjugation does not change the underlying real bundle of \(L^m\).
\end{proof}

The criterion is explicit in the weights.  With
\(\mu=[L_{\R}]-2\), equation \eqref{eq:real-line-Rm} gives
\[
 \gamma_{a,b}
 =\sum_{i=1}^{n+1}R_{|a_i|}(\mu)
  +\sum_{j=1}^{n+1}R_{|b_j|}(\mu).
\]
Thus the problem reduces to the relations in \(\widetilde{J}(\CP^n)\).  These can
be converted into finitely many congruence conditions on the weights,
but the formulas become unwieldy in higher dimensions; we therefore
retain \eqref{eq:general-J-homotopy} as the general statement and only
write out the low-dimensional consequences below.

For \(n=2,3\), put
\[
 \lambda(a,b)=\sum_i(a_i^2+b_i^2)=\Pi_1(a,b).
\]
In these two cases \(\widetilde{KO}(\CP^n)\cong\Z\), detected by the
coefficient of \(p_1\), and \cite[Lemma 1.4]{Adams1966} gives
\[
 \widetilde{J}(\CP^n)\cong\Z/24,
\]
with the stable \(J\)-homomorphism given by reduction modulo \(24\).  Rank \(2n+2\) is already in the stable range over
\(\CP^n\).

\begin{theorem}\label{thm:low-dimensional-homotopy}
Let \(n=2\) or \(3\), and let \(M_{a,b}\) and \(M_{a',b'}\) be the
corresponding orbit manifolds.  Then
\[
 M_{a,b}\simeq M_{a',b'}
 \quad\Longleftrightarrow\quad
 \lambda(a,b)\equiv\lambda(a',b')\pmod{24}.
\]
\end{theorem}

\begin{proof}
By \eqref{eq:eta-pontryagin},
$
 p_1(\eta_{a,b})=(\lambda(a,b)-(n+1))x^2.$
Hence congruence of the two \(\lambda\)-invariants modulo \(24\) is
exactly equality of the corresponding stable \(J\)-classes.  By
\cref{thm:sphere-realization}, the two orbit manifolds are PL
homeomorphic to the associated sphere bundles, which are then
fiber-homotopy equivalent.  This proves sufficiency.

Conversely, $p_1 \pmod{24}$ is a homotopy invariant for smooth manifolds. By $p_1(M_{a,b})=\lambda(a,b)x^2$, we finish the proof.
\end{proof}

For \(n=3\), formula \eqref{eq:eta-pontryagin} becomes
\[
 p_1(\eta_{a,b})=(\lambda(a,b)-4)x^2.
\]
In this dimension the sphere-bundle model can be chosen explicitly.

\begin{corollary}\label{cor:13-sphere-model}
For every \(13\)-dimensional orbit manifold \(M_{a,b}\), there are
integers \(q_1,q_2,q_3,q_4\) such that
\[
 q_1^2+q_2^2+q_3^2+q_4^2=\lambda(a,b)-4
\]
and
\[
 M_{a,b}\cong_{\PL}
 S\bigl((L^{q_1}\oplus L^{q_2}\oplus L^{q_3}\oplus L^{q_4})_{\R}\bigr).
\]
Equivalently, \(M_{a,b}\) is $\PL$ homeomorphic to the linear quotient
\(M_{\mathbf1,q}\), where \(\mathbf1=(1,1,1,1)\) and
\(q=(q_1,q_2,q_3,q_4)\).
\end{corollary}

\begin{proof}
The freeness criterion gives $\lambda(a,b)\ge4$: if one weight is zero,
all weights on the other sphere have absolute value one; if no weight is
zero, the assertion is immediate.  Lagrange's theorem supplies $q_i$ with
$\sum_iq_i^2=\lambda(a,b)-4$.  The action with weights $\mathbf1$ and $q$
is free, and the Hopf quotient on the first sphere identifies its total
space with
\[
 S\bigl((L^{q_1}\oplus L^{q_2}\oplus L^{q_3}\oplus L^{q_4})_{\R}\bigr).
\]
By \cref{prop:tangent}, its first Pontryagin coefficient is
$4+\sum_iq_i^2=\lambda(a,b)$.  Thus the claimed PL classification theorem,
\cref{thm:main}, applied to these \emph{two orbit manifolds}, gives the result.
\end{proof}

\section{Homogeneous Einstein and Sasakian geometry}\label{sec:Einstein}

We now turn to geometric consequences of the smooth classification.  For
coprime nonzero integers \(k,l\), let
\[
M_{k,l}^{(n)}
=
M_{(k,\ldots,k),(l,\ldots,l)}
=
\bigl(S^{2n+1}\times S^{2n+1}\bigr)/S^1,
\]
where each weight occurs \(n+1\) times.  The product Hopf fibration identifies
this quotient with the total space of a principal circle bundle $\xi$ over
\(\CP^n\times\CP^n\).  Up to interchanging the base factors and complex
conjugation, its Euler class is
\[
 e(\xi)=k\alpha_1+l\alpha_2,
\]
where \(\alpha_i\) is the positive generator of the second cohomology of the
\(i\)th factor.  We use this convention below and restrict to \(0<k\le l\).

\subsection{Homogeneous Einstein metrics}\label{subsec:Einstein}

Wang and Ziller \cite{WangZiller1990} constructed homogeneous Einstein
metrics on these circle bundles.  We first compute their Einstein constants
after volume normalization.

\begin{lemma}\label{lem:WZ-normalized-constant}
Let \(n\ge1\), let \(0<k\le l\) be coprime integers, and let
\(g_{k,l}\) be the Wang--Ziller Einstein metric on \(M_{k,l}^{(n)}\),
scaled so that its circle fibers have length \(2\pi\).  Put
\[
N=k^2+l^2.
\]
There is a unique parameter
\[
1\le r<\frac{n+1}{n}
\]
determined by \(k/l\) such that, with \(s=r+r^{-1}\), the Einstein
constant \(\widehat\rho_{k,l}\) of the corresponding unit-volume metric
satisfies
\begin{equation}\label{eq:normalized-Einstein}
\widehat\rho_{k,l}^{\,4n+1}
=
\frac{A_n}{N}
(s+2)^{2n+1}\bigl(2n+1-ns\bigr),
\end{equation}
where \(A_n>0\) depends only on \(n\).  Moreover, for fixed \(N\), distinct
unordered pairs \(\{k,l\}\) give distinct values of
\(\widehat\rho_{k,l}\).
\end{lemma}

\begin{proof}
Let \(g_1,g_2\) be the Fubini--Study metrics on the two copies of \(\CP^n\),
with K\"ahler forms normalized by
\[
[\omega_i]=2\pi\alpha_i,
\qquad
\Ric(g_i)=(n+1)g_i.
\]
Write \(\pi:M_{k,l}^{(n)}\to\CP^n\times\CP^n\) for the circle bundle and choose the
invariant connection \(\theta\) with
\[
\theta(e)=1,
\qquad
d\theta=\pi^*(k\omega_1+l\omega_2),
\]
where \(e\) generates the circle action with period \(2\pi\).  Suppressing
pullbacks, the metric has the form
\[
g_{k,l}=\theta^2+Xg_1+Yg_2,
\qquad X,Y>0,
\]
where \(\theta^2=\theta\otimes\theta\) denotes a symmetric $(0,2)$-tensor, which represents the metric contribution in the vertical directions tangent to the circle fibers of the principal $S^1$-bundle.
The Einstein constant $\rho$ of $g_{k,l}$ satisfies \cite[(1.11)--(1.12)]{WangZiller1990}
\begin{equation}\label{eq:WZ-specialized}
\begin{split}
\rho
&=\frac n2\left(\frac{k^2}{X^2}+\frac{l^2}{Y^2}\right),\\
\rho
&=\frac{n+1}{X}-\frac{k^2}{2X^2},\\
\rho
&=\frac{n+1}{Y}-\frac{l^2}{2Y^2}.
\end{split}
\end{equation}
Set \(r=X/Y\).  Substituting the last two equations into the first yields
\[
(2n+1)\rho X=n(n+1)(1+r).
\]
Substitution back gives
\begin{equation}\label{eq:WZ-X-rho}
X=\frac{(2n+1)k^2}{2(n+1)(n+1-nr)},
\qquad
Y=\frac Xr,
\qquad
\rho=\frac{n(n+1)(1+r)}{(2n+1)X}.
\end{equation}
The first formula may be rewritten as
\[
\frac{k^2}{X}
=
\frac{2(n+1)}{2n+1}(n+1-nr),
\]
while the third Einstein equation gives
\[
\frac{l^2r^2}{X}
=
\frac{2(n+1)}{2n+1}\bigl((n+1)r-n\bigr).
\]
Positivity implies
\begin{equation}\label{Interval}
\frac n{n+1}<r<\frac{n+1}{n}.
\end{equation}
Dividing the two identities yields
\begin{equation}\label{eq:k-l-r}
\frac{k^2}{l^2}
=
F_n(r)
:=
\frac{r^2(n+1-nr)}{(n+1)r-n}.
\end{equation}
Equivalently,
\begin{equation}\label{eq:WZ-r-cubic}
nl^2r^3-(n+1)l^2r^2+(n+1)k^2r-nk^2=0.
\end{equation}
On the interval \eqref{Interval},
\[
F_n'(r)
=
-\frac{r\bigl(2n(n+1)(r-1)^2+(2n-1)r\bigr)}
        {\bigl((n+1)r-n\bigr)^2}<0.
\]
Moreover,
\[
F_n(1)=1,
\qquad
\lim_{r\to (n+1)/n^-}F_n(r)=0.
\]
Hence
\[
F_n:\left[1,\frac{n+1}{n}\right)\longrightarrow(0,1]
\]
is a bijection, and \(0<k\le l\) determines a unique \(r\).

Put \(d=4n+1\).  The fixed fiber length and the base-volume normalization
give
\[
\Vol(g_{k,l})=B_n(XY)^n,
\qquad
B_n=\frac{(2\pi)^{2n+1}}{(n!)^2}.
\]
Thus the unit-volume metric has Einstein constant
\[
\widehat\rho_{k,l}=\rho\,\Vol(g_{k,l})^{2/d}.
\]
Using \eqref{eq:WZ-X-rho}, one obtains
\[
\widehat\rho_{k,l}^{\,d}
=
A_n(1+r)^d\frac{n+1-nr}{k^2}\,r^{-2n},
\]
where
\[
A_n
=
B_n^2\frac{2(n+1)}{2n+1}
\left(\frac{n(n+1)}{2n+1}\right)^{4n+1}.
\]
For \(s=r+r^{-1}\), equation \eqref{eq:k-l-r} gives
\[
\frac{n+1-nr}{k^2}
=
\frac{1+r}{Nr}\bigl(2n+1-ns\bigr).
\]
Together with \((1+r)^2/r=s+2\), this proves
\eqref{eq:normalized-Einstein}.

Finally, put
\[
\Phi_n(s)=(s+2)^{2n+1}(2n+1-ns).
\]
Then
\[
\Phi_n'(s)
=
(s+2)^{2n}
\bigl(1-2n-2n(n+1)(s-2)\bigr)<0
\qquad(s\ge2).
\]
For fixed \(N\), the normalized Einstein constant therefore determines
\(s\), hence \(r\), hence \(k^2/l^2\), and finally the unordered pair
\(\{k,l\}\).
\end{proof}

Let \(\mathcal E_h(M)\) denote the moduli space of unit-volume homogeneous
Einstein metrics on \(M\) with the quotient
\(C^\infty\)-topology; the transitive Lie group is not fixed in advance.

\begin{corollary}\label{cor:unbounded-Einstein-components}
Let \(p_1,\ldots,p_m\) be distinct primes congruent to \(1\pmod4\), and put
\[
N_m=p_1\cdots p_m.
\]
Then:
\begin{enumerate}[label=\textup{(\roman*)}]
\item every nine-dimensional orbit manifold with
\(\lambda=3N_m\) has
\[
\#\pi_0\mathcal E_h(M)\ge2^{m-1};
\]
\item every thirteen-dimensional orbit manifold with
\(\lambda=4N_m\) has
\[
\#\pi_0\mathcal E_h(M)\ge2^{m-1}.
\]
\end{enumerate}
In particular, the number of connected components of
\(\mathcal E_h(M)\) is unbounded in both the non-spin nine-dimensional
family and the spin thirteen-dimensional family.
\end{corollary}

\begin{proof}
The sum-of-two-squares formula gives exactly \(2^{m-1}\) primitive
representations
\[
N_m=k^2+l^2,
\qquad 0<k<l.
\]
All are primitive because \(N_m\) is square-free, and \(N_m\) is odd.

For \(n=2\),
\[
p_1\bigl(M_{k,l}^{(2)}\bigr)=3N_mx^2,
\qquad
w_2\bigl(M_{k,l}^{(2)}\bigr)=\bar x\ne0.
\]
Thus the non-spin nine-dimensional classification identifies every orbit
manifold with \(\lambda=3N_m\) with each of the spaces
\(M_{k,l}^{(2)}\).  For \(n=3\), similarly,
\[
p_1\bigl(M_{k,l}^{(3)}\bigr)=4N_mx^2,
\qquad
w_2\bigl(M_{k,l}^{(3)}\bigr)=0,
\]
and \cref{thm:spin4} identifies every orbit manifold with
\(\lambda=4N_m\) with each \(M_{k,l}^{(3)}\).

Transport the Wang--Ziller metrics by these diffeomorphisms and normalize
volumes to one.  By \cref{lem:WZ-normalized-constant}, the
\(2^{m-1}\) representations give \(2^{m-1}\) distinct positive Einstein
constants.  By \cite[proof of Theorem~1.8]{BohmWangZiller2004}, for every
\(\varepsilon>0\) only finitely many Einstein constants occur among
unit-volume homogeneous Einstein metrics on a fixed homogeneous compact connected manifold with
scalar curvature at least \(\varepsilon\).
 Suppose that two such metrics with distinct positive Einstein constants \(\hat{\rho}_1<\hat{\rho}_2\) lie in the same connected component of \(\mathcal E_h(M)\).  Scalar curvature defines a continuous function on \(\mathcal E_h(M)\), i.e., for $g\in \mathcal E_h(M)$,
$$\mathrm{Scal}(g)=d\hat{\rho}(g)$$ 
where \(d=\dim M\), $\hat{\rho}(g)$ is the Einstein constants.
 The image of this connected component contains the entire interval \([d\hat{\rho}_1,d\hat{\rho}_2]\). This would yield infinitely many distinct positive Einstein constants with scalar curvature at least \(d\hat{\rho}_1>0\), contradicting the finiteness result above. Hence, metrics with distinct positive Einstein constants belong to distinct connected components of \(\mathcal E_h(M)\).
  This proves
the stated lower bounds.
\end{proof}

\subsection{Regular Sasakian and contact structures}\label{subsec:Sasakian}

The same repeated-weight presentations carry another natural geometric
structure.  Assume \(k,l>0\).  The class
\[
k\alpha_1+l\alpha_2\in H^2(\CP^n\times\CP^n;\Z)
\]
is represented by the integral K\"ahler form
\(k\omega_1+l\omega_2\).  The Boothby--Wang construction therefore equips
\(M_{k,l}^{(n)}\) with a regular Sasakian structure
\cite{BoothbyWang1958,BoyerGalicki2008}.  Let
\(\mathcal D_{k,l}\) denote its contact distribution.

\begin{proposition}\label{prop:contact-c1}
For coprime positive integers \(k,l\), choose a generator
\(x\in H^2(M_{k,l}^{(n)};\Z)\) so that
$
\pi^*\alpha_1=l x,$
$\pi^*\alpha_2=-k x.$
Then
\begin{equation}\label{eq:contact-c1}
c_1(\mathcal D_{k,l})=(n+1)(l-k)x.
\end{equation}
\end{proposition}

\begin{proof}
For a regular Sasakian Boothby--Wang fibration, the complex structure
on the contact distribution is the horizontal lift of the complex
structure on the K\"ahler base.  Equivalently,
\[
(\mathcal D_{k,l},J)
\cong
\pi^*T(\CP^n\times\CP^n)
\]
as complex vector bundles; see
\cite[Section~2.1]{BoyerTonnesenFriedman2012}.
Consequently,
\[
c_1(\mathcal D_{k,l})
=
\pi^*c_1\bigl(T(\CP^n\times\CP^n)\bigr).
\]
This pullback formula is also written explicitly in
\cite[Section~4.3, equation~(24)]{BoyerTonnesenFriedman2012}.

The Gysin sequence gives
\[
H^2(M_{k,l}^{(n)};\Z)
\cong
\frac{\Z\alpha_1\oplus\Z\alpha_2}
     {\langle k\alpha_1+l\alpha_2\rangle}
\cong\Z,
\]
and the displayed choice of \(x\) realizes this quotient. Therefore, 
\[
\pi^*\alpha_1=l x,
\qquad
\pi^*\alpha_2=-k x.
\]
Since
\[
c_1\bigl(T(\CP^n\times\CP^n)\bigr)
=(n+1)(\alpha_1+\alpha_2),
\]
formula \eqref{eq:contact-c1} follows.
\end{proof}

\begin{corollary}\label{cor:unbounded-Sasakian}
With \(N_m=p_1\cdots p_m\) as in
\cref{cor:unbounded-Einstein-components}, every nine-dimensional orbit
manifold with \(\lambda=3N_m\), and every thirteen-dimensional orbit
manifold with \(\lambda=4N_m\), carries at least \(2^{m-1}\) regular
Sasakian structures whose underlying cooriented contact structures are
pairwise non-contactomorphic.
\end{corollary}

\begin{proof}
For each primitive representation
\[
N_m=k^2+l^2,
\qquad 0<k<l,
\]
transport the Boothby--Wang structure on \(M_{k,l}^{(n)}\) to the fixed
smooth manifold using the diffeomorphisms supplied by the smooth
classification.  Suppose two of the resulting contact structures were
contactomorphic.  Since \(H^2(M;\Z)\cong\Z\), the underlying diffeomorphism
acts on \(H^2\) by multiplication by \(\pm1\).  Hence
\cref{prop:contact-c1} would give
\[
|l-k|=|l'-k'|.
\]
Because \(0<k<l\) and \(0<k'<l'\), write this common difference as \(d>0\).
Then
\[
N_m=k^2+(k+d)^2=k'^2+(k'+d)^2.
\]
The right-hand side is strictly increasing in the nonnegative integer
variable, so \(k=k'\) and \(l=l'\).  Thus distinct representations give
non-contactomorphic contact structures.  There are \(2^{m-1}\) such
representations.
\end{proof}


\section*{Declaration on the Use of AI Tools}

During the preparation of this manuscript, the authors used ChatGPT (OpenAI) as an auxiliary tool for improving the English language, including grammar correction, refinement of mathematical exposition, and enhancement of clarity and readability. ChatGPT was also used to assist in developing two geometric applications, concerning the moduli spaces of unit-volume homogeneous Einstein metrics and regular Sasakian structures with pairwise non-contactomorphic underlying contact structures. In addition, AI-assisted computations were employed in the proof of Proposition \ref{prop:delta-weights}.

The research problems, central ideas, and main conceptual framework of this work were conceived by the authors. All AI-generated suggestions and computational outputs were critically examined, and all mathematical arguments and results were checked and verified by the authors, who take full responsibility for the correctness and integrity of the work.

\end{document}